\documentclass[12pt,a4paper]{amsart}
\usepackage{amssymb}
\usepackage{amsfonts}
\usepackage{amsthm}
\usepackage{amsmath}
\usepackage{amscd}
\usepackage[latin2]{inputenc}
\usepackage{t1enc}
\usepackage[mathscr]{eucal}
\usepackage{indentfirst}
\usepackage{graphicx}
\usepackage{graphics}
\numberwithin{equation}{section}
\usepackage[margin=2.9cm]{geometry}
\usepackage{epstopdf} 
\usepackage{xcolor}
 
\usepackage{verbatim}
 \usepackage{enumitem}
\usepackage[numbers,sort&compress]{natbib}

\theoremstyle{plain}
\newtheorem{theorem}{Theorem}[section]
\newtheorem{lemma}[theorem]{Lemma}

\newtheorem{proposition}[theorem]{Proposition}

 \theoremstyle{definition}

\newtheorem{remark}[theorem]{Remark}
\newtheorem{?}[theorem]{Problem}

\newcommand{\R}{\mathbb{R}}
\newcommand{\II}{{\rm{II}}}
\newcommand{\D}{\mathbb{D}}

\newcommand{\A}{\mathcal A}

\begin{document}
	
	\title[Eigenvalue Variations]{Eigenvalue Variations of the Neumann Laplace Operator Due to Perturbed Boundary Conditions}

	\author[M. Nursultanov]{Medet Nursultanov}
	\address {Department of Mathematics and Statistics, University of Helsinki}
	\email{medet.nursultanov@gmail.com}

	\author[W. Trad] {William Trad}
	\address {School of Mathematics and Statistics, University of Sydney}
	\email{w.trad@maths.usyd.edu.au}
	
	\author[J. Tzou] {Justin Tzou}
	\address {School of Mathematical and Physical Science, Macquarie University}
	\email{tzou.justin@gmail.com}
	
	\author[L. Tzou] {Leo Tzou}
	\address {University of Amsterdam}
	\email{leo.tzou@gmail.com}



	\begin{abstract}
		This work considers the Neumann eigenvalue problem for the weighted Laplacian on a Riemannian manifold $(M,g,\partial M)$ under the singular perturbation. This perturbation involves the imposition of vanishing Dirichlet boundary conditions on a small portion of the boundary. We derive a sharp asymptotic of the perturbed eigenvalues, as the Dirichlet part shrinks to a point $x^*\in \partial M$, in terms of the spectral parameters of the unperturbed system. This asymptotic demonstrates the impact of the geometric properties of the manifold at a specific point $x^*$. Furthermore, it becomes evident that the shape of the Dirichlet region holds significance as it impacts the first terms of the asymptotic. A crucial part of this work is the construction of the singularity structure of the restricted Neumann Green's function which may be of independent interest. We employ a fusion of layer potential techniques and pseudo-differential operators during this work.
	\end{abstract}
 
	\maketitle

	\tableofcontents
	
	\section{Introduction}
	Comprehending the disturbances within physical fields caused by inhomogeneities in a known environment is crucial for various purposes. It helps to understand the robustness of the body's behaviour under small perturbations of its constituent material; see \cite{AmmariKang2004, BonnetierDapognyVogelius} for more such applications. Mathematically, this task involves performing an asymptotic analysis of the solution of the partial differential equation, when defining domain or properties of the material are slightly perturbed. Many examples of this general question have been studied, including investigations into the conductivity equation \cite{CapdeboscqVogelius,AmmariSeo,Cedio-FengyaMoskowVogelius} and other areas, such as linearized elasticity \cite{AmmariKangNakamuraTanuma, BerettaBonnetierFranciniMazzucato}, Maxwell equations \cite{AmmariVogeliusVolkov,Griesmaier}, and cellular biology \cite{holcman2015stochastic,BressloffNewby}. 
	
	Here, we investigate the eigenvalues of the weighted Laplace operator under mixed Dirichlet-Neumann boundary conditions when the Dirichlet region disappears. More precisely, we formulate the problem as follows. Let $(M,g)$ be a compact, connected, orientable Riemannian manifold with smooth non-empty boundary $\partial M$. Consider the eigenvalue problem
	\begin{equation}\label{eig_val_prob}
		-\Delta_g u - g(F,\nabla_g u) = \lambda u, \qquad \left.u\right|_{\Gamma_\varepsilon} = 0,\qquad \left.\partial_\nu u\right|_{\partial M\setminus\Gamma_\varepsilon} = 0,
	\end{equation}
	where $\Delta_g$ is the negative Laplace-Beltrami operator, $\nabla_g$ is the gradient, $\nu$ is the outward pointing normal vector field, $F$ is a force field, and $\Gamma_\varepsilon\subset \partial M$ is a connected piece of boundary of size $\varepsilon>0$. We denote the corresponding operator by $-\Delta_{Mix,\varepsilon}^F$. The objective is to derive an asymptotic of the eigenvalues $\{\lambda_{j,\varepsilon}\}_{j\in \mathbb{N}}$ of $-\Delta_{Mix,\varepsilon}^F$ as $\varepsilon$ tends to zero, that is, as $\Gamma_\varepsilon$ shrinks in a suitable sense that will be specified later. We will do this in terms of the spectral parameters of the unperturbed operator, which is the weighted Neumann Laplacian, denoted by $-\Delta_N^F$.
	
	The problem at hand is closely related to the "narrow escape problem," which has gained significant attention in recent years due to its relevance in cellular biology \cite{holcman2015stochastic,BressloffNewby}. In this scenario, $M$ denotes a cavity with a reflecting boundary, except for a small absorbing window $\Gamma_\varepsilon$. The particles in $M$ are modelled as Brownian motions that exit only through the region $\Gamma_\varepsilon$. The mean first-passage time, which represents the expected duration a particle will wander before escaping, is a crucial metric in this context. The narrow escape problem concerns the asymptotic behaviour of the mean first-passage time as the size of the window $\Gamma_\varepsilon$ tends towards zero.
	
	\subsection{Previous results} The investigation of the behaviour of eigenvalues of elliptic boundary value problems under the singular perturbation of boundary conditions has a long history, see for instance \cite{AbatangeloFelliLena,ammari,Borisov2003Izvest,BorisovSibMatZh,Cherdantsev,FelliNorisOgnibene,FelliNorisOgnibene2022,GadylshinDiffUrav,Gadylshin1993,Gadylshin1998AlgAnal,GadylshinRepevskiiShishkinaJMS,nursultanovTradTzouTzou,Planida2005,EswarathasanKolokolnikov}. Detailed analysis of the two-dimensional planar domain has been performed in \cite{Gadylshin1993}, where the author provided the full asymptotic expansion of the perturbed eigenvalues. Moreover, a complete pointwise expansion of the perturbed eigenfunctions is provided as well.
	
	In \cite{Cherdantsev}, the perturbed eigenvalues in a three-dimensional Euclidean domain with a smooth boundary were studied. The author derived the asymptotic behaviour of the perturbed eigenvalues up to an unspecified $o(\varepsilon)$ term:
	\begin{equation*}
		\lambda_{j,\varepsilon} = \lambda_j + 4\pi |u_j(x^*)|^2 c\varepsilon + o(\varepsilon),
	\end{equation*}
	where $\{\lambda_j\}_{j\in \mathbb{N}}$ and $\{u_j\}_{j\in \mathbb{N}}$ are unperturbed eigenvalues and the corresponding normalized eigenfunctions, $c$ is the constant depending on geometry of $\Gamma_\varepsilon$.
	
	The problem of perturbed eigenvalues in a domain with a Lipschitz boundary in the Euclidean space was examined by the authors in \cite{FelliNorisOgnibene}. They obtained the asymptotic behaviour of the perturbed eigenvalues, expressed in terms of the unperturbed eigenvalues and the relative Sobolev $u_j$-capacity of $\Gamma_\varepsilon$:
	\begin{equation*}
		\lambda_{j,\varepsilon} = \lambda_j + \mathrm{Cap}_{M}(\Gamma_\varepsilon, u_j) + o\left(\mathrm{Cap}_{M}(\Gamma_\varepsilon, u_j)\right).
	\end{equation*}
	
	The case of a two-dimensional planar domain in the presence of a force field was considered in \cite{ammari}. The authors used layer potential techniques to derive the asymptotic expansion
    \begin{equation*}
		\lambda_{j,\varepsilon} = \lambda_j + \pi|u_j(x^*)|^2e^{\phi(x^*)}\log\varepsilon+O(|\log\varepsilon|^{-2}),
    \end{equation*}
    where $\phi$ is a potential, that is $F = \nabla \phi$.

    Additionally, we refer to the studies \cite{WardKeller,KolokolnikovTitcombeWard,CoombsStraubeWard,CheviakovWard,PaquinLefebvreIyaniwuraWard}, which employed the 
    method of matched asymptotic expansions.
	\subsection{Main results}
	Despite the large number of works on this topic, there are still many questions regarding more general geometries. The present work is devoted to answering this question. In our setting, the Dirichlet region, $\Gamma_{\varepsilon,a}$, is considered to be a small geodesic ellipse centred at $x^*\in\partial M$, with eccentricity $\sqrt{1 - a^2}$ and size $\varepsilon \rightarrow 0$ (to be made precise later). Moreover, the force field is given by a smooth up to the boundary potential $\phi$, that is $F=\nabla_g\phi$. Our main objective is to investigate how the geometric characteristics of the manifold $M$ at a specific point $x^*$ impact the asymptotic expansion. We use a combination of the methods used in \cite{ammari} and \cite{NursultanovTzouTzou}: layer potential technique and microlocal analysis.
	
	In our current analysis, it is crucial to understand the singular structure of the Neumann Green function, which is represented by the following equation:
	\begin{equation*}
		\begin{cases}
			\Delta_gG_M^\omega(x,y)-\text{div}_g(F(y)G_M^\omega(x,y))+\omega^2G_M^\omega(x,y) = -\delta_x(y), \\
			\left.\partial_\nu G_M^\omega(x,y) - g_y(F(y),\nu)G_M^\omega(x,y)\right\vert_{y\in \partial M} = 0.
		\end{cases}
	\end{equation*}
	where $\omega^2$ is an element of the resolvent set of the eigenvalue problem \eqref{eig_val_prob}. Let us consider a formal restriction of $G_M^\omega$ to the boundary $\partial M$. We denote this by the symbol $G_{\partial M}^\omega$. For an exact definition of this restriction, refer to Section \ref{Green_function_sec}. We obtain the singularities structure of $G_{\partial M}^\omega$ near the diagonal and in the neighbourhood of an eigenvalue of \eqref{eig_val_prob}, when $w^2$ approaches it. To state it, let us set the necessary notions. For $x$, $y\in \partial M$, let $H(x)$ denote the mean curvature of the boundary at $x$, $d_g(x,y)$ the geodesic distance given by metric $g$, $d_h(x,y)$ the geodesic distance given by induced metric $h$ on the boundary $\partial M$, and
	\begin{equation*}
		\mathrm{II}_x(V): = \mathrm{II}(V,V), \qquad V\in T_x\partial M,
	\end{equation*}
	the scalar second fundamental form (see pages 235 and 381 of \cite{lee} for definitions).
	\begin{proposition}\label{sing_stract}
		Let $(M,g,\partial M)$ be a compact connected orientable Riemannian manifold of dimension three with a non-empty smooth boundary. Let $\lambda_j$ be a simple eigenvalue of $-\Delta_N^F$ and $V_j$ be a neighbourhood of $\lambda_j$ which does not contain any other eigenvalue of $-\Delta_N^F$. Then there exists a neighbourhood of 
		$$\mathrm{Diag}:=\{(x,x)\in \partial M\times \partial M\}$$ 
		where the singularity structure of $G_{\partial M}^\omega$ given by:
		\begin{align}\label{sing_stract_formula}
		\nonumber	G^\omega_{\partial M}(x,y)&=\frac{1}{2\pi}d_g(x,y)^{-1} -\frac{H(x)}{4\pi}\log d_h(x,y)+\frac{g_x(F,\nu)}{4\pi}\log d_h(x,y)\\ \nonumber&+\frac{1}{16\pi}\left(\mathrm{II}_{x}\left(\frac{\exp_{x}^{-1}(y)}{|\exp_{x}^{-1}(y)|_h}\right)-\mathrm{II}_{x}\left(\frac{\star\exp_{x}^{-1}(y)}{|\exp_{x}^{-1}(y)|_h}\right)\right) \\
			& +\frac{1}{4\pi}h_{x}\left(F^{||}(x),\frac{\exp^{-1}_{x}(y)}{|\exp^{-1}_{x}(y)|_h}\right)\\
		\nonumber	& + \frac{u_j(x)u_j(y)}{\lambda_j - \omega^2}e^{\phi(y)}  + R_{\partial M}^\omega(x,y),
		\end{align}
		where $\omega^2 \in V_j\setminus \{\lambda_j\}$, $R^\omega_{\partial M}(x,y) \in C^{0,\alpha}(\partial M\times \partial M)$ for $\alpha\in (0,1)$, $F^{||}$ denotes the tangential component of $F$ and $\star$ denotes the Hodge star operator. 
	\end{proposition}
    In deriving the singularity structure, we employed a standard pseudo-differential parametrix construction as described in \cite{NursultanovTzouTzou, NursultanovTradTzou} by observing that $G_{\partial M}^\omega$ is an approximate inverse to a Dirichlet-to-Neumann map. To determine the singularity with respect to $\omega$ near the eigenvalues of $-\Delta_N^F$, we have used an approach similar to that of \cite{HHHmonog}. Using the aforementioned proposition, we derive an asymptotic expression for the eigenvalues $\lambda_{j,\varepsilon}$ as $\varepsilon\rightarrow 0$. For the sake of clarity, we begin by presenting the main result for $\Gamma_{\varepsilon,a}$ being a geodesic ball, that is $a=1$:
	
	\begin{theorem}\label{main_for_ball}
		Let $(M,g,\partial M)$ be a compact connected orientable Riemannian manifold of dimension three with a non-empty smooth boundary. Fix $x^*\in \partial M$ and let $\Gamma_{\varepsilon}$ be the boundary geodesic ball centred at $x^*$ of geodesic radius $\varepsilon>0$. Assume that $F = \nabla_g \phi$ for a potential $\phi$ smooth up to the boundary. Let $\{\lambda_{j,\varepsilon}\}_{j\in \mathbb N}$ be the eigenvalues of $-\Delta_{Mix,\varepsilon}^F$. If $\lambda_j$ is a simple eigenvalue of $-\Delta_N^F$ and $u_j$ is the corresponding eigenfunction normalized in $L^2(M,e^\phi d\mu_g)$ (weighted $L^2$ space with a weight $e^\phi$), then
		\begin{equation*}
			\lambda_{j,\varepsilon} - \lambda_j = A\varepsilon + B\varepsilon^2 \log\varepsilon + C\varepsilon^2 + O(\varepsilon^3\log^2\varepsilon),
		\end{equation*}
		where
		\begin{align*}
			&A = 4 |u_j(x^*)|^2e^{\phi(x^*)},\\[1em]
			&B = 4 \pi |u_j(x^*)|^2e^{\phi(x^*)} (H(x^*)-\partial_{\nu}\phi(x^*)),\\[1em]
			&C = |u_j(x^*)|^2 e^{\phi(x^*)} \left( \frac{8\log 2 - 6}{\pi}(H(x^*)-\partial_{\nu}\phi(x^*)) - 16R_{\partial M}^{\lambda_j}(x^*,x^*) \right).
		\end{align*}
		Here, $H$ is the mean curvature of the boundary, $R_{\partial M}^{\lambda_j}(x^*,x^*)$ is the evaluation at $(x,y)=(x^*,x^*)$ of the kernel $R_{\partial M}^{\lambda_j}(x,y)$ in Proposition \ref{sing_stract}.
	\end{theorem}
	Theorem \ref{main_for_ball} does not realize the full power of Proposition \ref{sing_stract} as it does not see the inhomogeneity of the local geometry at $x^*$, only the mean curvature shows up. This limitation arises from the fact that Dirichlet regions are specifically geodesic balls. However, by considering geodesic ellipses instead of geodesic balls, we observe that the inclusion of the second fundamental form term in Proposition \ref{sing_stract} contributes to an asymptotic term, which is the difference in principal curvatures. The ellipse, we consider, is defined as follows. Let $E_1(x^*)$, $E_2(x^*) \in T_{x^*}\partial M$ be the unit eigenvectors of the shape operator at $x^*$ corresponding respectively to principal curvatures $\kappa_1(x^*)$ and $\kappa_2(x^*)$. For $a\in (0,1]$ fixed, we set
	\begin{align}\label{ellipse_intro}
		\Gamma_{\varepsilon,a}:=\{\mathrm{exp}_{x^*;h}(\varepsilon t_1 E_1+\varepsilon t_2 E_2)\mid t_1^2+a^{-2}t_2^2\leq 1\}.
	\end{align}
	Now, we are ready to state the main result:
	\begin{theorem}\label{main_for_ellipse}
		Let $(M,g,\partial M)$ be a compact connected orientable Riemannian manifold of dimension three with a non-empty smooth boundary. Fix $x^*\in \partial M$ and let $\Gamma_{\varepsilon,a}$ be the boundary geodesic ellipse given by \eqref{ellipse_intro}. Assume that $F = \nabla_g \phi$ for a potential $\phi$ smooth up to the boundary. Let $\{\lambda_{j,\varepsilon}\}_{j\in \mathbb N}$ be the eigenvalues of $-\Delta_{Mix,\varepsilon}^F$. If $\lambda_j$ is a simple eigenvalue of $-\Delta_N^F$ and $u_j$ is the corresponding eigenfunction normalized in $L^2(M,e^\phi d\mu_g)$, then
		\begin{equation*}
			\lambda_{j,\varepsilon} - \lambda_j = A\varepsilon + B\varepsilon^2 \log\varepsilon + (C_1 + C_2 + C_3)\varepsilon^2 + O(\varepsilon^3\log^2\varepsilon).
		\end{equation*}
		Here, the constants are defined as follows
		\begin{align*}
			& K_a = \frac{\pi}{2} \int_0^{2\pi} \frac{a}{(a^2\cos^2\theta + \sin^2\theta)^{1/2}}d\theta\\
			& A = \frac{4\pi^2a}{K_a} |u_j(x^*)|^2e^{\phi(x^*)},\\
			& B = \frac{4\pi^3 a^2 |u_j(x^*)|^2e^{\phi(x^*)}}{K_a^2} (H(x^*)-\partial_{\nu}\phi(x^*))
		\end{align*}
		and
		\begin{multline*}
			C_1 = \frac{a^2\pi(H(x^*)-\partial_{\nu}\phi(x^*))|u_j(x^*)|^2          e^{\phi(x^*)}}{K_a^2}\times\\
			\times\int_{\mathbb{D}} \frac{1}{ (1-|s'|^2)^{1/2}}  \int_\D  \frac{\log\left((t_1 - s_1)^2 + a^2 (t_2-s_2)^2 \right)^{1/2}}{(1-|t'|^2)^{1/2}} dt' ds',
		\end{multline*}
		\begin{multline*}
			C_2 = \frac{a^2\pi\left(\kappa_1(x^*) -                         \kappa_2(x^*)\right)|u_j(x^*)|^2 e^{\phi(x^*)} }{4K_a^2}\times\\
			\times \int_{\mathbb{D}} \frac{1}{ (1-|s'|^2)^{1/2}} \int_\D \frac{(t_1 - s_1)^2 - a^2 (t_2 - s_2)^2}{(t_1 - s_1)^2 + a^2 (t_2 - s_2)^2} \frac{1}{ (1-|t'|^2)^{1/2}} dt' ds'
		\end{multline*}
		
		\begin{flalign*}
			\;\:\;\: C_3 = \frac{16 \pi^4 a^2R_{\partial M}^{\lambda_j}(x^*,x^*) |u_j(x^*)|^2 e^{\phi(x^*)} }{K_a^2}, &&
		\end{flalign*}   
		where $R_{\partial M}^{\lambda_j}(x^*,x^*)$ is the evaluation at $(x,y)=(x^*,x^*)$ of the kernel $R_{\partial M}^{\lambda_j}(x,y)$ in Proposition \ref{sing_stract}.
	\end{theorem}
	From this result, we can conclude that the shape of the Dirichlet region is important. The eccentricity of the ellipse affects the main term of the asymptotic. Moreover, we observe a difference in the principal curvatures, which is not visible in case $a=1$. 
	
	\subsection{Outline of the paper} In Section \ref{notations}, we initiate the exposition by introducing the necessary notations and providing a more precise formulation of the problem at hand. Section \ref{Green_function_sec} is devoted to the computation of the singular structure of Green's function. Moving on to Section \ref{proof_main_result}, we employ variational principles to derive additional bounds for the perturbed eigenvalues. Ultimately, in this section, we establish the proof of the main theorem.

	\section{Preliminaries}\label{notations}
	In this section, we introduce basic notations and formulate the problem. Throughout this work, we use $(M, g)$ to denote a compact connected orientable Riemannian manifold of dimension three with a non-empty smooth boundary. The corresponding volume form and geodesic distance are denoted by $d\mu_g(\cdot)$ and $d_g(\cdot,\cdot)$, respectively. Let $\iota_{\partial M}: \partial M \hookrightarrow M$ be the trivial embedding of the boundary $\partial M$ into $M$. This allows us to define the boundary metric $h := \iota_{\partial M}^*g$ inherited by $g$. We similarly use $d\mu_h(\cdot)$ and $d_h(\cdot,\cdot)$ to denote respectively the volume form on the boundary and the geodesic distance on the boundary given by metric $h$. We denote the Laplace-Beltrami operator by $\Delta_g = -d^*d$.
	
	For $x\in \partial M$, let $E_1(x), E_2(x) \in T_{x}\partial M$ be the unit eigenvectors of the shape operator at $x\in \partial M$ corresponding respectively to the principal curvatures $\kappa_1(x),\ \kappa_2(x)$. We will drop the dependence in $x$ from our notation when there is no ambiguity. We choose $E_1$ and $E_2$ such that $E_1^\flat\wedge E_2^\flat\wedge\nu^\flat$ is a positive multiple of the volume form $ d\mu_g$ (see p.26 of \cite{lee} for the ``musical isomorphism'' notation of $^\flat$ and $^\sharp$). Here we use $\nu$ to denote the outward-pointing normal vector field. By $H(x)$, we denote the mean curvature of $\partial M$ at $x$. We also set 
	$$\II_x(V) := \II_x(V,V),\ \ V\in T_x\partial M,$$
	to be the scalar second fundamental form. Note that, in defining $\II$ and the shape operator, we will follow the standard literature in geometry (e.g. \cite{lee}) and use the inward-pointing normal so that the sphere embedded in $\R^3$ would have positive mean curvature in our convention.
	
	In this article, we will often use boundary normal coordinates. Therefore, we briefly recall its construction. For a fixed $x^*\in \partial M$, we will denote by $B_h(\rho;x^*) \subset \partial M$ the geodesic disk of radius $\rho>0$ (with respect to the metric $h$) centred at $x^*$ and $\D_\rho$ to be the Euclidean disk in $\R^2$ of radius $\rho$ centred at the origin. In what follows $\rho$ will always be smaller than the injectivity radius of $(\partial M, h)$. Letting $t = (t_1,t_2,t_3) \in \R^3$, we will construct a coordinate system $x(t; x^*)$ by the following procedure:
	
	Write $t\in \R^3$ near the origin as $t  = (t',t_3)$ for $t' = (t_1,t_2)\in \D_\rho$. Define first 
	$$x((t',0); x^*) := {\rm {exp}}_{x^*;h} (t_1 E_1+ t_2 E_2),$$
	where ${\rm{exp}}_{x^*;h} (V)$ denotes the time $1$ map of $h$-geodesics with initial point $x^*$ and initial velocity $V\in T_{x^*}\partial M$. The coordinate $t'\in \D_\rho \mapsto x((t',0); x^*)$ is then an $h$-geodesic coordinate system for a neighborhood of $x^*$ on the boundary surface $\partial M$. We can then construct a coordinate system for a neighbourhood of $x^*\in M$ by considering $g$-geodesic rays $\gamma_{x^*,-\nu}:[0,\rho)\rightarrow M$ emanating from points in $\partial M$ orthogonal to $\partial M$. In particular, we can then smoothly extend $t'$ to $U$ by setting $t'$ to be constant functions along $\gamma_{x^*,-\nu}$. If we then define $t_3$ to be the unit speed parameter of $\gamma_{x^*,-\nu}$, then $(t_1,t_2,t_3)$ form coordinates for $M$ in some neighborhood of $x^*\in M$. As a consequence, $t_3$ is a boundary-defining function, that is $t_3>0$ away from $\partial M$ and $t_3=0$ on $\partial M$. We will call these local coordinates, \textit{boundary normal coordinates}. For convenience we will write $x(t'; x^*)$ in place of $x((t',0);x^*)$. Readers wishing to know more about boundary normal coordinates can refer to \cite{lee} for a brief recollection of the basic properties we use here for detailed construction.  

	We will also use the rescaled version of this coordinate system. For $\varepsilon >0$ sufficiently small we define the (rescaled) $h$-geodesic coordinate by the following map
	\begin{equation}\label{res coord}
		x^{\varepsilon}(\cdot ; x^*) : t' = (t_1, t_2) \in \D \mapsto x(\varepsilon t'; x^*) \in B_h(\varepsilon;x^*),
	\end{equation}
	where $ \mathbb{D}$ is the unit disk in $\R^2$. Given the boundary normal co-ordinate construction, we define the geodesic ellipse $\Gamma_{\varepsilon,a}$ as the following subset of $\partial M$ 
	\begin{align}\label{ellipse}
		\Gamma_{\varepsilon,a}:=\{\mathrm{exp}_{x^*;h}(\varepsilon t_1 E_1+\varepsilon t_2 E_2)\mid t_1^2+a^{-2}t_2^2\leq 1\}.
	\end{align}
	\subsection{Formulation of the problem}
	Now we are ready to state the problem. Let us consider the operator
	\begin{equation*}
		u \rightarrow \Delta_gu + g(F,\nabla_gu),
	\end{equation*}
	where $F$ is a force field, which is given by $F=\nabla_g\phi$ for a smooth up to the boundary potential $\phi$. We can re-write this operator in the following way
	\begin{equation*}
		\Delta_{g}^F\cdot := \Delta_g\cdot + g(F,\nabla_g\cdot) = \frac{1}{e^\phi}\mathrm{div}_g(e^\phi \nabla_{g} \cdot).
	\end{equation*}
	According to  \cite{GrigoryanSaloff-Coste}, the operator $\Delta_{g}^F$ is called a weighted Laplacian and the pair $(M,e^\phi d\mu_g)$ is called a weighted manifold. Note that $e^\phi$ is bounded and strictly positive on $M$. Therefore, $L^2(M) = L^2(M,e^\phi d\mu_g)$ as sets with equivalent norms. We also note that the operator $\Delta_{g}^{F}$ with initial domain $C_0^\infty(M)$ is essentially self-adjoint in $L^2(M,e^\phi d\mu_g)$ and non-positive definite. 
	We want to study the operator $\Delta_{g}^F$ with Dirichlet and Neumann boundary conditions on $\Gamma_{\varepsilon, a}$ and $\partial M\setminus \Gamma_{\varepsilon, a}$, respectively.  This operator can be defined via quadratic form as follows. Consider the quadratic form
	\begin{equation*}
		a_{\varepsilon}(u,v) :=  \int_{M} e^{\phi(z)} g(\nabla_g u(z), \nabla_g u(z)) d\mu_g(z),
	\end{equation*}
	with the domain 
	\begin{equation*}
		\mathrm{D}(a_\varepsilon) := \{u \in H^1(M): \mathrm{supp}(u\arrowvert_{\partial M} \subset \partial M\setminus \Gamma_{\varepsilon, a} \}.
	\end{equation*}
	Since $e^{\phi}$ is strictly positive and bounded on $M$, this quadratic form is non-negative, closed, symmetric and densely defined in $L^2(M,e^\phi d\mu_g)$, and hence, generates the self-adjoint non-negative operator; see Theorem 2.6 in \cite[Ch. 6.2]{Kato}. We denote this operator by $-\Delta^F_{Mix,\varepsilon}$ and call it a weighted Laplace operator corresponding to the aforementioned mixed boundary conditions.
	
	Since $H^1(M)$ is compactly embedded in $L^2(M,e^\phi d\mu_g)$, the spectrum of $-\Delta^F_{Mix,\varepsilon}$ is discrete and consists of the eigenvalues with finite multiplicity accumulating at infinity. We denote them, taking into account their multiplicities, as follows
	\begin{equation*}
		0 \leq \lambda_{1,\varepsilon}\leq \lambda_{2,\varepsilon} \leq  \cdots < \infty.
	\end{equation*}
	The corresponding normalized, in the $L^2(M,e^\phi d\mu_g)$ sense, eigenfunctions are denoted by $\{u_{j,\varepsilon}\}_{j\in \mathbb{N}}$.
	
	We also consider the operator $-\Delta_g^F$ with a Neumann boundary condition, which is generated by the quadratic form
	\begin{equation*}
		a_N(u,u) :=  \int_{M} e^{\phi(z)} g(\nabla_g u(z), \nabla_g u(z)) d\mu_g(z), \qquad \text{with } \mathrm{D}(a_N) = H^1(M).
	\end{equation*}
	We denote this operator by $-\Delta_N^F$. By the same arguments above, we conclude that the spectrum $\mathrm{spec}(-\Delta_N^F)$ is discrete and consists of eigenvalues with finite multiplicity accumulating at infinity. We  donate them by 
	\begin{equation*}
		0=\lambda_{1} \leq \lambda_{2}\leq \cdots < \infty.
	\end{equation*}
	By $\{u_{j}\}_{j\in \mathbb{N}}$, we denote the corresponding normalized, in the $L^2(M,e^\phi d\mu_g)$ sense, eigenfunctions. We aim to derive an asymptotic expansion for $\lambda_{j,\varepsilon}$ as $\varepsilon\rightarrow 0$ in terms of $\lambda_j$ and $u_j$.

\section{Neumann Greens Function}\label{Green_function_sec}
In this section, we consider the Greens function $G_M^\omega$ defined as the solution (in the distributional sense) to the following boundary value problem
\begin{align*}
		\begin{cases}
			\Delta_gG_M^\omega(x,y)-\text{div}_g(F(y)G_M^\omega(x,y))+\omega^2G_M^\omega(x,y) = -\delta_x(y), \\
			\left.\partial_\nu G_M^\omega(x,y) - g_y(F(y),\nu)G_M^\omega(x,y)\right\vert_{y\in \partial M} = 0,
		\end{cases}
\end{align*}
where $\omega^2$ is a parameter belonging to the resolvent set of the operator $-\Delta_N^F$. We seek the singularity structure of $\partial M$-restriction of $G_M^\omega$ near the diagonal. A more precise definition of this restriction will be given later. Our main aim is to obtain Proposition \ref{sing_stract}. 
\subsection{Singularity in the spectral parameter}
The first step in our analysis is to express $G_M^\omega$ as a decomposition of Neumann eigenfunctions. The following result is similar to a result of \cite{HHHmonog}. We modify it here for our setting
	\begin{proposition}\label{sing on omega}
		Let $\{\lambda_j\}_{j\in \mathbb{N}}$ be the eigenvalues of $-\Delta_N^F$ and $\{u_j\}_{j\in \mathbb{N}}$  be the corresponding  $L^2(M,e^\phi d\mu_g)$-orthonormalized eigenfunctions. Then, for $x\neq y$ and $\omega^2\in\mathbb{C}\setminus \mathrm{spec}(-\Delta_N^F)$, it follows that
		\begin{equation*}
			G_M^\omega(x,y) = \sum_{j=1}^\infty \frac{u_j(x)u_j(y)}{\lambda_j-\omega^2}e^{\phi(y)}.
		\end{equation*}
	\end{proposition}
\begin{proof}
		Since $(M,g,\partial M)$ is assumed to be compact, it follows that $e^\phi$ is bounded and strictly positive on $M$. Therefore, $L^2(M) = L^2(M,e^\phi d\mu_g)$ with equivalent norms. Since $G_M^\omega(x,\cdot)  \in L^2(M,e^\phi d\mu_g)$, for any $x\in M$, we can express 
		\begin{equation*}
			G_M^\omega(x,y) = \sum_{j=1}^\infty v_j(x) u_j(y) e^{\phi(y)}.
		\end{equation*}
		Since $e^\phi >0$ is bounded, it follows that $f\in L^2(M)$ if and only if $e^\phi f\in L^2(M)$. Thus, for fixed $x\in M$, the above expression for $G_M^\omega(x,y)$ is unique. Green's identity in conjunction with the divergence theorem as well as the boundary condition on $u_j$ yields the following calculation
		\begin{align*}
			u_j(x) &= \int_M \left(-\Delta_g G_M^\omega(x,y) + \text{div}_y(F(y)G_M^\omega(x,y)) - \omega^2G_M^\omega(x,y)\right)u_j(y) d\mu_g(y),\\
			&= (\lambda_j-\omega^2) \sum_{k=1}^\infty \int_M v_k(x) u_k(y) u_j(y) e^{\phi(y)} d\mu_g(y), \\ 
			&= (\lambda_j - \omega^2) v_j(x) \int_M |u_j(y)|^2 e^{\phi(y)} d\mu_g(y). 
		\end{align*}
		Recall that $u_j$ are $L^2(M,e^\phi d\mu_g)$-orthonormalized so that $v_j(x) = u_j(x)/(\lambda_j-\omega^2)$, which implies that 
		\begin{align*}
			G_M^\omega(x,y) = \sum_{j = 1}^\infty \frac{u_j(x)u_j(y)}{\lambda_j-\omega^2}e^{\phi(y)},
		\end{align*}
		as required. 
	\end{proof}
Let $\lambda_j$ be a simple eigenvalue of $-\Delta_N^F$ and $V_j$ be an open bounded neighborhood of $\lambda_j$ in $\mathbb{C}$ such that $V_j\cap \mathrm{spec}(-\Delta_N^F) = \{\lambda_j\}$. Within Section \ref{Green_function_sec}, we are interested in deriving explicit asymptotics for the \textit{trace} of $G^\omega_M$. To this end, we write $G^\omega_M$ as 
	\begin{align*}
		G^\omega_M(x,y) &= \sum_{k\neq j} \frac{u_k(x)u_k(y)}{\lambda_k-\omega^2}e^{\phi(y)} + \frac{u_j(x)u_j(y)}{\lambda_j-\omega^2}e^{\phi(y)},\\
		&:= N^\omega_M(x,y) + \frac{u_j(x)u_j(y)}{\lambda_j-\omega^2}e^{\phi(y)},
	\end{align*}
	for $\omega^2\in V_j$. From here we can then define $G^\omega_{\partial M}$ and $N^\omega_{\partial M}$ as the boundary restrictions of $G^\omega_M$ and $N^\omega_M$ as Schwartz kernels of the trace of the integral operators $G^\omega_M$ and $N^\omega_M$ respectively. That is, for $f\in C^\infty(\partial M)$, we have 
	\begin{align*}
		&G^\omega_{\partial M}:f\mapsto \left.\left(\int_{\partial M} G^\omega_M(x,y)f(y)d\mu_h(y)\right)\right\vert_{x\in \partial M},\\
		&N^\omega_{\partial M}:f\mapsto \left.\left(\int_{\partial M} N^\omega_M(x,y)f(y)d\mu_h(y)\right)\right\vert_{x\in \partial M}.
	\end{align*}
	It can also be seen from the perspective of microlocalization that the above restrictions to $\partial M \times \partial M$ are well-defined since $G^\omega_M$ is the Schwartz kernel of a pseudo-differential operator with $\mathrm{WF}(G^\omega_M)\cap N^*(\partial M\times \partial M)=\emptyset$ and the difference between $G^\omega_M$ and $N^\omega_M$ is a $C^\infty(\partial M\times \partial M)$ term. Thus, we write the trace of $G^\omega_M$ as 
	\begin{align}
		G^\omega_{\partial M}(x,y) = N^\omega_{\partial M}(x,y) +\frac{u_j(x)u_j(y)}{\lambda_j-\omega^2}e^{\phi(y)}.\label{Restricted greens}
	\end{align}
 Our choice of $N^\omega_M$ suggests that any terms which depend on $\omega^2$ in $N^\omega_{\partial M}$ are negligible, since, by definition $\omega^2\in V_j\setminus \{\lambda_j\}$ and $V_j$ is judiciously chosen such that there are no other Neumann eigenvalues in $V_j$. This implies that there is only one significant singularity in $\omega^2$ which is given by the second term of \eqref{Restricted greens}. 
\subsection{Singularities along the diagonal}
 When considering \eqref{Restricted greens}, it is apparent that $\frac{u_j(x)u_j(y)}{\lambda_j-\omega^2}e^{\phi(y)}$ is jointly smooth on $\partial M\times \partial M$, thus the only difficulty in deriving asymptotics for $G^\omega_{\partial M}$ for $x$ near $y$ lies in the derivation of $N^\omega_{\partial M}$. We will show that $N_{\partial M}^\omega$ is a left parametrix for a Dirichlet-to-Neumann map, which is associated to the following auxiliary Dirichlet boundary value problem
    \begin{equation}\label{Auxiliary Problem}
        \begin{cases}
            \Delta_g u_f +g(F,\nabla_g u_f)+\omega^2u_f = 0,\\
            u_f\arrowvert_{\partial M} = f\in C^\infty(\partial M).
        \end{cases}
    \end{equation}
    The Dirichlet-to-Neumann map is given by 
    \begin{align*}
        &\Lambda_{g,F}^\omega: H^{1/2}(\partial M)\ni f \mapsto \partial_\nu u_f \in H^{1/2}(\partial M)^*.
    \end{align*}
In order to construct $N^\omega_{\partial M}$, we will require a series of technical lemmas. The first of which was proven in \cite{LLSConformalLaplacian} and \cite{LeeUhlmann89}. We offer a sketch of the proof for our special case under consideration.
\begin{lemma}\label{DtN Map lemma}
    The Dirichlet-to-Neumann map $\Lambda^\omega_{g,F}$ is an elliptic pseudo-differential operator of order $1$. In addition, the first two terms of the symbol of $\sigma(\Lambda^\omega_{g,F})(t,\xi')$ are 
    \begin{align*}
			\sigma_1(\Lambda_{g,F}^\omega) &= -\sqrt{\widetilde{q_2}}, \\
			\sigma_0(\Lambda_{g,F}^\omega) &= \frac{1}{2\sqrt{\widetilde{q_2}}}(\nabla_{\xi'}\sqrt{\widetilde{q_2}}\cdot D_{t'}\sqrt{\widetilde{q_2}}-\widetilde{q_1}-\partial_{t_3}\sqrt{\widetilde{q_2}} +\widetilde{E}\sqrt{\widetilde{q_2}}),
    \end{align*}
    where $\widetilde{E}$, $\widetilde{q_1}$ and $\widetilde{q_2}$ are given by 
    \begin{align*}
			\widetilde{E}(t) &:= -\frac{1}{2}\sum_{\alpha,\beta} h^{\alpha\beta}(t)\partial_{t_3}h_{\alpha\beta}(t)-F^3(t),\\
			\widetilde{q_1}(t,\xi')&:=-i\sum_{\alpha,\beta} \left(\frac{1}{2}h^{\alpha\beta}(t)\partial_{t_\alpha}\log \delta(t) +\partial_{t_\alpha} h^{\alpha\beta}(t) - F^\alpha(t)h^\beta_\alpha(t)\right)\xi_\beta,\\ 
			\widetilde{q_2}(t,\xi')&:=\sum_{\alpha,\beta} h^{\alpha\beta}(t)\xi_\alpha\xi_\beta,
    \end{align*}
    Here $\alpha,\beta\in \{1,2\}$,  $F=F^1(t)\partial_{t_1}+F^2(t)\partial_{t_2}+F^3(t)\partial_{t_3}$ and $\delta(t)=\mathrm{det}(g_{\alpha\beta})$.
\end{lemma}
\begin{proof}
		Within our choice of co-ordinates, we begin with the following decomposition
		\begin{align}
			-\Delta_g - g(F,\nabla_g\cdot)-\omega^2 = D_{t_3}^2 + i \widetilde{E}(t) D_{t_3} + \widetilde{Q^\omega}(t,D_{t'}), \label{Helmholtz op}
		\end{align}
		where $\widetilde{Q^\omega}(t,D_{t'})$ is 
		\begin{align*}
			\widetilde{Q^\omega}(t,D_{t'}) &:= \sum_{\alpha,\beta}h^{\alpha\beta}(t)D_{t_\alpha}D_{t_\beta} \\ 
			&-i\sum_{\alpha,\beta}\left(\frac{1}{2}h^{\alpha\beta}(t)\partial_{t_\alpha}\log\delta(t)+\partial_{t_\alpha}h^{\alpha\beta}(t)-F^\alpha(t)h^\beta_\alpha(t)\right)D_{t_\beta} -\omega^2.
		\end{align*}
		It should be noted that the total symbol $\sigma(\widetilde{Q^\omega})(t,\xi')$ of $\widetilde{Q^\omega}(t,D_{t'})$ is given by 
		\begin{align*}
			\sigma(\widetilde{Q^\omega})(t,\xi') &= \widetilde{q_1}(t,\xi') + \widetilde{q_2}(t,\xi')-\omega^2.
		\end{align*}
		It follows that we can construct a first order, classical, pseudo-differential operator $A^\omega_F(t,D_{t'})$ such that 
		\begin{align}
			-\Delta_g - g(F,\nabla_g\cdot)-\omega^2  &= (D_{t_3}+i\widetilde{E}(t)-iA^\omega_F(t,D_{t'}))(D_{t_3}+iA^\omega_F(t,D_{t'})),\label{Helmholtz decomposition}
		\end{align}
       by equating \eqref{Helmholtz op} and \eqref{Helmholtz decomposition}. This yields the following equation, up to a smoothing operator
		\begin{equation}\label{eq111111}
			A^\omega_F(t,D_{t'})^2 +i[D_{t_3},A^\omega_F(t,D_{t'})] - \widetilde{Q^\omega}(t,D_{t'}) - \widetilde{E}(t)A^\omega_F(t,D_{t'})  = 0.
		\end{equation}
		The standard pseudo-differential calculus allows us to write \eqref{eq111111} equivalently in terms of symbols associated with the relevant operators (ones which involve some action on functions defined over $\partial M$) as follows, up to a smoothing symbol  
		\begin{align*}
			\sum_{\gamma} \frac{1}{\gamma!} \partial_{\xi'}^\gamma \sigma(A^\omega_F) D_t^\gamma\sigma(A^\omega_F)-\partial_{t_3} \sigma(A^\omega_F) - \sigma(\widetilde{Q^\omega}) - \widetilde{E}(t) \sigma(A^\omega_F) = 0.
		\end{align*}
		where $\gamma\in \mathbb{N}^n$ denotes some multi-index. Collecting homogeneous terms of degree $2$ and then $1$ yields the first two terms of the Borel expansion for the symbol of the pseudo-differential operator $A^\omega_F(t,D_{t'})$. These terms are 
    \begin{align*}
			\sigma_1(A^\omega_F) &= -\sqrt{\widetilde{q_2}}, \\
			\sigma_0(A^\omega_F) &= \frac{1}{2\sqrt{\widetilde{q_2}}}(\nabla_{\xi'}\sqrt{\widetilde{q_2}}\cdot D_{t'}\sqrt{\widetilde{q_2}}-\widetilde{q_1}-\partial_{t_3}\sqrt{\widetilde{q_2}} +\widetilde{E}\sqrt{\widetilde{q_2}}).
    \end{align*}
    It should be noted that $\sigma_1(A^\omega_F)=0$ \textit{only if} $\xi'=0$ (which corresponds to the zero section). Thus, it follows that $A^\omega_F(t,D_{t'})$ is an elliptic operator. Furthermore, by construction, $\sigma_1(A^\omega_F)$ and $\sigma_0(A^\omega_F)$ are homogeneous symbols, as are the residual symbols $\sigma_{j}(A^\omega_F)$ for $j\leq-1$. Therefore $A^\omega_F(t,D_{t'})$ is a first-order, elliptic, classical pseudo-differential operator. Within the following section of the proof, it is shown that $A^\omega_F(t,D_{t'})$ coincides with $\Lambda_{g,F}^\omega$ up to a smoothing operator. This is done by first considering the region $\{0\leq t_3\leq T\}$. The authors in \cite{LeeUhlmann89}, exploited \eqref{Helmholtz decomposition} in order to write \eqref{Auxiliary Problem} as a system of forward and backward heat equations:
    \begin{align*}
			(D_{t_3}+iA^\omega_F)u_f &= v, \;\ \left.u_f\right\vert_{t_3=0}=f, \\
			(D_{t_3}+i\widetilde{E}-iA^\omega_F)v &= w\in C^\infty([0,T];\mathcal{D}'(\mathbb{R}^2)),
    \end{align*}
    where $u_f,v\in C^\infty([0,T];\mathcal{D}'(\mathbb{R}^{2}))$. Since $\sigma_1(A^\omega_F) < 0$ for $\xi'\neq 0$, the following heat equation is well-posed 
		\begin{align*}
			\partial_{t_3}v+A^\omega_Fv-\widetilde{E}v = -iw,
		\end{align*}
		and thus $v\in C^\infty([0,T]\times \mathbb{R}^{2})$. Thus, restricting the forward heat equation to $\{t_3 = 0\}$ implies that, up to a smoothing operator
		\begin{align*}
			\Lambda_{g,F}^\omega f = \left.\partial_{t_3} u\right\vert_{t_3=0} = \left.A_F^\omega(t,D_{t'}) u \right\vert_{t_3=0}. 
		\end{align*} 
    Consequently, we have that  
    $$\sigma(\Lambda_{g,F}^\omega)(t,\xi') = \sigma(A^\omega_F)(t,\xi').$$
\end{proof}

In addition to establishing the above lemma, we require a lemma that we can use to link $N^\omega_{\partial M}$ and $\Lambda^\omega_{g,F}$. The following lemma elucidates said link as it shows that $N^\omega_{\partial M}$ is a left, elliptic pseudo-differential parametrix of $\Lambda^\omega_{g,F}$. Once again, we consider the elliptic Dirichlet boundary value problem \eqref{Auxiliary Problem}. 
    \begin{lemma}\label{Parametrix lemma}
        The Dirichlet-to-Neumann map $\Lambda^\omega_{g,F}$ and $N^\omega_{\partial M}$ satisfy the following operator equation
		\begin{align}\label{Operator Form Eq}
			I = N^\omega_{\partial M} \Lambda^\omega_{g,F} +\Psi^{-\infty},
		\end{align}
        where $\Psi^{-\infty}$ denotes the class of smoothing operators. In particular,  $N^\omega_{\partial M} \in \Psi^{-1}_{cl}$ is an elliptic pseudo-differential operator. 
    \end{lemma}
	\begin{proof}
		We prove the above lemma by integrating by parts \eqref{Auxiliary Problem} against the Neumann Greens function $G^\omega_M(x,y)$
    \begin{align*}
        -u_f(x) &= \int_M G^\omega_M(x,z) \Delta_g u_f d\mu_g(z) + \int_{\partial M} \left(u_f(z)\partial_\nu G^\omega_M(x,z) - G^\omega_M(x,z)\partial_\nu u_f(z)\right) d\mu_h(z)\\
        &- \int_M u_f(z) \text{div}_g(F(z)G^\omega_M(x,z)) d\mu_g(z)+\omega^2\int_M u_f(z) G^\omega_M(x,z) d\mu_g(z).
    \end{align*}
    The divergence theorem yields
    \begin{align*}
        -u_f(x) &= \int_M G^\omega_M(x,z) \Delta_g u_f d\mu_g(z) + \int_{\partial M} u_f(z)\partial_\nu G^\omega_M(x,z) - G^\omega_M(x,z)\partial_\nu u_f(z) d\mu_h(z) \\
			&+\int_M g_z(F(z),\nabla_g u_f) G^\omega_M(x,z) d\mu_g(z) - \int_{\partial M} u_f(z)G^\omega_M(x,z) F(z) \cdot \nu d\mu_h(z)  \\
			&+ \omega^2 \int_M u_f(z) G^\omega_M(x,z)d\mu_g(z).
    \end{align*}
    Employing the prescribed Neumann boundary conditions on $G^\omega_M$, we conclude 
		\begin{align*}
			u_f(x) &= \int_{\partial M} G^\omega_M(x,z) \partial_\nu u_f(z) d\mu_h(z).
		\end{align*}
	Furthermore, restricting $x\in \partial M$ and invoking the prescribed Dirichlet boundary condition from \eqref{Auxiliary Problem} using our definition involving the trace \eqref{Restricted greens}, we have that 
		\begin{align*}
			f(x) &= \int_{\partial M} G^\omega_{\partial M}(x,z) \Lambda^\omega_{g,F}f(z) d\mu_h(y), \\ 
			&= \int_{\partial M} N^\omega_{\partial M}(x,z) \Lambda^\omega_{g,F}f(z) d\mu_h(z) + \int_{\partial M} e^{\phi(z)}\frac{u_j(x)u_j(z)}{\lambda_j-\omega^2} \Lambda_{g,F}^\omega f(z) d\mu_h(z).
		\end{align*}
    Since the Neumann eigenfunctions are smooth, the rightmost integral in the above expression consists of a smooth Schwartz kernel and thus gives rise to a smoothing operator. That is, we have 
		\begin{align*}
			I = N^\omega_{\partial M}\Lambda^\omega_{g,F} + \Psi^{-\infty}
		\end{align*}
    Finally, since $\Lambda_{g,F}^\omega \in \Psi^1_{cl}$ is an elliptic operator, we conclude that $N^\omega_{\partial M}\in \Psi^{-1}_{cl}$. 
\end{proof}
Using Lemmas \ref{DtN Map lemma} and \ref{Parametrix lemma} we can prove the following theorem by iteratively determining the terms in the Borel summation associated with the parametrix $N_{\partial M}^\omega$ modulo smoothing terms. 

 \begin{proposition}\label{sing_structure_for_N}
		Let $x,y\in \partial M$ such that $x\neq y$ and $\omega^2\in \mathbb{C}\setminus \mathrm{spec}(-\Delta_N^F)$. In an open neighbourhood of $\mathrm{Diag}:=\{(x,x)\in \partial M\times \partial M\}$, we have that
		\begin{align*}
			N^\omega_{\partial M}(x,y) &= \frac{1}{2\pi}d_g(x,y)^{-1} -\frac{H(x)}{4\pi}\log d_h(x,y)+\frac{g_x(F,\nu)}{4\pi}\log d_h(x,y)\\ &+\frac{1}{16\pi}\left(\mathrm{II}_{x}\left(\frac{\exp_{x}^{-1}(y)}{|\exp_{x}^{-1}(y)|_h}\right)-\mathrm{II}_{x}\left(\frac{\star\exp_{x}^{-1}(y)}{|\exp_{x}^{-1}(y)|_h}\right)\right) \\
			& +\frac{1}{4\pi}h_{x}\left(F^{||}(x),\frac{\exp^{-1}_{x}(y)}{|\exp^{-1}_{x}(y)|_h}\right) + R^\omega_{\partial M}(x,y),
		\end{align*}
		where $R^\omega_{\partial M}(x,y) \in C^{0,\alpha}(\partial M\times \partial M)$ for $\alpha\in (0,1)$, and $F^{||}$ denotes the tangential component of $F$ and $\star$ denotes the Hodge star operator. 
	\end{proposition}

Within the above expression for $N^\omega_{\partial M}$, we obtain a clear image as to the structure of the singularity in \(x,y\in \partial M\) near the diagonal (where \(x=y\)) as well as the singularity in \(\omega^2\in \mathbb{C}\setminus \mathrm{spec}(-\Delta_N^F)\) for $\omega^2$ near $\lambda_j$. For the sake of clarity, we will include an outline for the proof of Theorem \ref{sing_structure_for_N}. Since we have already determined the nature of the leading order singularity in $\omega$, all that is left is to reveal the nature of the singularity of the leading order terms for \(x\) near \(y\) on \(\partial M\). 
\begin{proof}[Proof of Proposition \ref{sing_structure_for_N}]
	There are infinitely many additional terms in the asymptotic series of $N^\omega_{\partial M}$, these are formed via an iterative argument on the level of symbols. However, for our purposes, we only need the first two elements of the kernel expansion and thus, we only need the first two symbols in the Borel expansion of $\sigma(N^\omega_{\partial M})$. Upon deriving $\sigma_{-1}(N^\omega_{\partial M})$ and $\sigma_{-2}(N^\omega_{\partial M})$ An expression for the asymptotic series of the Schwartz kernel is given by the fourier transform of $\sigma_{-1}(N^\omega_{\partial M})+\sigma_{-2}(N^\omega_{\partial M})$. To begin this iterative process, we view the operator equation \eqref{Operator Form Eq} on the level of symbols. 
		\begin{equation}
			1 = \sigma(N^\omega_{\partial M})\#\sigma(\Lambda_{g,F}^\omega)(x,\xi') + S^{-\infty}, \label{SymbolCalcEq}
		\end{equation}
		where $\#$ denotes the standard composition of symbols in correspondence to the composition of pseudo-differential operators. Furthermore, we write $\sigma(N^\omega_{\partial M})(t,\xi')$ in terms of the following asymptotic series, whose existence is guaranteed by Borel's lemma
		\begin{align*}
			\sigma(N^\omega_{\partial M})(x,\xi') \sim \sum_{j\geq 1} \sigma_{-j}(N^\omega_{\partial M})(t,\xi'), \;\ \sigma_{-j}(N^\omega_{\partial M})(t,\xi') \in S^{-j}_{1,0}.
		\end{align*}
		Equation \eqref{SymbolCalcEq} becomes the following in accordance to the formula for the $\#$-product
		\begin{align*}
			1 = \sum_{\gamma}\frac{1}{\gamma!}\partial_{\xi'}^\gamma \sigma(N^\omega_{\partial M}) D_t^\gamma \sigma(\Lambda_{g,F}^\omega) + S^{-\infty}.
		\end{align*} 
		where $\gamma\in\mathbb{N}^n$ denotes a multi-index. In the first iteration, we have that for a smooth function $\chi$ and $R\in \mathbb{R}$ with $\chi(\xi') = 0$ for $|\xi'| \leq R$ and $\chi(\xi')=1$ for $|\xi'|\geq 2R$
		\begin{align*}
			1 = \sigma_{-1}(N^\omega_{\partial M})(t,\xi')\sigma_1(\Lambda_{g,F}^\omega)(t,\xi') + S^{-1}_{1,0} \implies \sigma_{-1}(N^\omega_{\partial M})(t,\xi') = \frac{\chi(\xi')}{\sigma_1(\Lambda_{g,F}^\omega)(t,\xi')}.
		\end{align*}
		We can further iterate for the second term by forming the following equation
		\begin{align*}
			1 &= \sigma_{-1}(N^\omega_{\partial M})(t,\xi')\sigma_1(\Lambda_{g,F}^\omega)(t,\xi')+\sigma_{-1}(N^\omega_{\partial M})(t,\xi')\sigma_0(\Lambda_{g,F}^\omega)(t,\xi'),\\&+\sigma_{-2}(N^\omega_{\partial M})(t,\xi')\sigma_1(\Lambda_{g,F}^\omega)(t,\xi') +\nabla_{\xi'} \sigma_{-1}(N^\omega_{\partial M}) \cdot D_{t'} \sigma_1(\Lambda_{g,F}^\omega) + S^{-2}_{1,0}.
		\end{align*}
		We now equate terms of symbol order $-1$ to obtain the following equation
		\begin{align*}
			0 &= \sigma_{-1}(N^\omega_{\partial M})(t,\xi')\sigma_0(\Lambda_{g,F}^\omega)(t,\xi')+\sigma_{-2}(N^\omega_{\partial M})(t,\xi')\sigma_1(\Lambda_{g,F}^\omega)(t,\xi')\\&+\nabla_{\xi'} \sigma_{-1}(N^\omega_{\partial M}) \cdot D_{t'} \sigma_1(\Lambda_{g,F}^\omega).
		\end{align*}
		So, we choose $\sigma_{-2}(N^\omega_{\partial M})(t,\xi')$ as follows
		\begin{align*}
			\sigma_{-2}(N^\omega_{\partial M})(t,\xi') = -\frac{\chi(\xi')}{\sigma_1(\Lambda_{g,F}^\omega)(t,\xi')}\left(p_{-1}(t,\xi')\sigma_0(\Lambda_{g,F}^\omega)(t,\xi')+\nabla_{\xi'} p_{-1} \cdot D_{t'} \sigma_1(\Lambda_{g,F}^\omega)\right).
		\end{align*}
		Thus, we have that the Schwartz kernel of the $\partial M$-restricted Greens function, when evaluated at the center of the $h$-geodesic disc $\Gamma_{\varepsilon,1}$, up to $\Psi^{-3}$, when written in local co-ordinates is given by 
		\begin{align*}
			\frac{1}{4\pi^2}\left(\int_{\mathbb{R}^{2}}e^{-i\xi'\cdot t'}\sigma_{-1}(N^\omega_{\partial M})(0,\xi')d\xi' + \int_{\mathbb{R}^{2}}e^{-i\xi'\cdot t'}\sigma_{-2}(N^\omega_{\partial M})(0,\xi') d\xi'\right).
		\end{align*}
		This results in the desired singular expansion for the boundary-restricted Greens function fixed at a central point $x^*$. It can then be extended via a series of estimates derived in \cite{NursultanovTzouTzou} to \(N^\omega_{\partial M}(x,y)\) for \(x\neq y\), but suitably close (See \cite{NursultanovTradTzou} for full calculation).
\end{proof}

Finally, we note that Proposition \ref{sing_stract} follows as a trivial consequence of Proposition \ref{sing_structure_for_N} and relation \eqref{Restricted greens}

\subsection{Schwartz kernel estimates} Within this section, we investigate $G_{\partial M}^\omega$ near $x^*$, in the local coordinates given by \eqref{res coord}. We introduce several integral operators related to the terms on the right-hand side of \eqref{sing_stract_formula}. First, we consider a weighted variant of the normal operator
\begin{equation}\label{La}
    L_a f = a \int_{\mathbb{D}} \frac{f(s')}{\left((t_1 - s_1)^2 + a^2(t_2 - s_2)^2\right)^{1/2}} ds'
\end{equation}
acting on functions on the disk $\D$. It is known that $L_a$ is a self-adjoint operator; see for instance Section 4 in \cite{NursultanovTzouTzou}. Moreover, by \cite{schuss2006ellipse}, it follows that 
\begin{equation}\label{La u = 1}
    L_a  \left({K_a}^{-1} {(1 - |t'|^2)^{-1/2}}\right) = 1, \;\ K_a = \frac{\pi}{2} \int_{0}^{2 \pi} \left( \cos^2 \theta + \frac{\sin^2 \theta}{a^2} \right)^{-1/2} d \theta. 
\end{equation}
	By (4.4) in \cite{NursultanovTzouTzou}, it was shown that $u(t')=K_a^{-1}(1-|t'|^2)^{-1/2}$ is the unique solution to $L_a u = 1$ in $H^{1/2}(\mathbb{D})^*$. Next, we introduce the following operators
	\begin{align*}
		R_{\log,a} f(t') &:= a \int_\D \log\left((t_1 - s_1)^2 + a^2 (t_2-s_2)^2 \right)^{1/2} f(s') ds', \\ 
		R_{\infty,a} f(t') &:= a \int_{\mathbb{D}} \frac{(t_1 - s_1)^2 - a^2 (t_2 - s_2)^2}{(t_1 - s_1)^2 + a^2 (t_2 - s_2)^2} f(s') ds',\\
		R_{F,a} f(t') &:= a \int_{\mathbb{D}} \frac{F^1(0)(t_1 - s_1) + a F^2(0) (t_2 - s_2)}{((t_1 - s_1)^2 + a^2 (t_2 - s_2)^2)^{1/2}} f(s') ds',\\
		R_{I,a}f(t') &:= a\int_{\mathbb{D}} f(s')ds'.
	\end{align*}
	\begin{remark}\label{rem}
		In \cite{NursultanovTzouTzou}, it was shown that the operators $R_{\log,a}$ and $R_{\infty,a}$ are bounded maps from $H^{1/2}(\D)^*$ to $H^{3/2}(\D)$. Repeating the arguments shows that this is also true for $R_{F,a}$. 
	\end{remark}
 Note that these lemmas are proved in \cite{NursultanovTzouTzou,NursultanovTradTzou}. We state them here for the convenience of the reader.
\begin{lemma}\label{l1}
		We have the following identity
		\begin{equation*}
			\int_{\Gamma_{\varepsilon,a}} d_g(x,y)^{-1}v(y)d\mu_h(y) = \varepsilon L_a\tilde v (t') + \varepsilon^3 \mathcal{A}_\varepsilon\tilde v(t'),
		\end{equation*}
		where $x=x^{\varepsilon}(t')$, $\tilde v(t') = v(x^{\varepsilon}(t'))$, for some $\A_\varepsilon : H^{1/2}(\D; ds')^* \to H^{1/2}(\D; ds')$ with operator norm bounded uniformly in $\varepsilon$.
\end{lemma}

From now on, we will denote by $\A_\varepsilon$ any operator which takes 
$$\A_\varepsilon: H^{1/2}(\D ;ds')^* \to H^{1/2}(\D; ds'),$$ 
whose operator norm is bounded uniformly in $\varepsilon$. 
	\begin{lemma}\label{l2}
		The following identity holds
		\begin{multline*}
			(H(x)-g_x(F,\nu))\int_{\Gamma_{\varepsilon,a}} \log d_h(x,y) v(y)d\mu_h(y)\\
			= \varepsilon^2\log\varepsilon(H(x^*)-\partial_{\nu}\phi(x^*)) R_I\tilde v(t') + \varepsilon^2 (H(x^*)-\partial_{\nu}\phi(x^*))R_{\log,a} \tilde v(t') + \varepsilon^3\log\varepsilon \mathcal{A}_\varepsilon \tilde v(t'),
		\end{multline*}
		where $x=x^{\varepsilon}(t')$ and $\tilde v(t') = v(x^{\varepsilon}(t'))$.
	\end{lemma}
\begin{lemma}
		The following identity holds
		\begin{multline*}
			\int_{\Gamma_{\varepsilon,a}}\left(\mathrm{II}_x\left(\frac{\exp_x^{-1}(y)}{|\exp_x^{-1}(y)|_h}\right)-\mathrm{II}_x\left(\frac{\star\exp_x^{-1}(y)}{|\exp_x^{-1}(y)|_h}\right)\right)v(y)d\mu_h(y)\\
			=\varepsilon^2(\kappa_1(x^*)-\kappa_2(x^*)) R_{\infty,a} \tilde v(t') + \varepsilon^3 \A_\varepsilon \tilde v(t'),
		\end{multline*}
		where $x=x^{\varepsilon}(t')$ and $\tilde v(t') = v(x^{\varepsilon}(t'))$. Recall that $\kappa_1(x^*)$ and $\kappa_2(x^*)$ denote the principle curvatures of the boundary $\partial M$ at $x^*$. 
	\end{lemma}

\begin{lemma}\label{l4}
		The following identity holds
		\begin{align*}
			\int_{\Gamma_{\varepsilon,a}} h_x\left(F^{\parallel}(x), \frac{\exp_{x; h}(y)}{|\exp_{x; h}(y)|_h}\right) v(y)d\mu_h(y) = \varepsilon ^2 R_{F,a}\tilde v(t') + \varepsilon^3\mathcal{A}_\varepsilon\tilde v(t'),
		\end{align*}
		where $x=x^{\varepsilon}(t')$ and $\tilde v(t') = v(x^{\varepsilon}(t'))$.
\end{lemma}
	
Finally, we need to know the behaviour of the final component on the right-hand side of equation \eqref{sing_stract_formula} as the parameter $\omega$ converges towards the spectrum of $-\Delta_N^F$. eigenvalue of $-\Delta_g^F$.

\begin{proposition}\label{tr}
		Let $\lambda_k$ be a simple eigenvalue of $-\Delta_g^F$ and let $V_k$ be its neighbourhood which is open, bounded, and does not contain any other eigenvalue of $-\Delta_g^F$. For $\lambda\in V_k$, let
		\begin{equation}
			R_{\lambda_k,\lambda}: C^\infty(\partial M) \mapsto \mathcal{D}'(\partial M)
		\end{equation}
		be the operator defined by the integral kernel 
		$$R^{\lambda_k}(x,y) - R^\lambda(x,y),$$
		then
		\begin{equation*}
			\left\| R_{\lambda_k,\lambda} \right\|_{H^{1/2}(\partial M)^* \mapsto H^{1/2}(\partial M)} = O(|\lambda_k - \lambda|).
		\end{equation*}
\end{proposition} 

\begin{remark}
    Due to Proposition \ref{sing on omega}, for any fixed $x$, $y\in \partial M$, $R_{\partial M}^{\lambda_k}(x,y)$ is well defined.
\end{remark}
	
\begin{proof}[Proof of Proposition \ref{tr}]
Throughout the proof, we write $x\lesssim y$ or $y\gtrsim x$ to mean that $x\leq Cy$, where $C >0$ is some constant.
The dependencies of $C$ will be clear from the context. By $x\approx y$ we mean that $x\lesssim y$ and $x\gtrsim y$.

    For $\psi\in H^{1/2}(\partial M)^*$, we have the following estimate
    \begin{equation*}
			\left\| \int_{\partial M} (R^{\lambda_k}(x,y) - R^\lambda(x,y))\psi(y)d\mu_h(y)\right\|_{H^{1/2}(\partial M)} \leq \|U\|_{H^1(M)},
    \end{equation*}
		where
    \begin{equation*}
			U(x) = \sum_{j\neq k} \frac{(\lambda_k - \lambda)u_j(x)\langle u_je^\phi,\psi\rangle}{(\lambda_j - \lambda_k)(\lambda_j - \lambda)}
    \end{equation*}
    and $\langle \cdot,\cdot \rangle$ denotes paring between $H^{1/2*}$ and $H^{1/2}$. Let us consider the following Neumann boundary value problem
    \begin{equation}\label{hol_1}
			\begin{cases}
				(-\Delta_g^F + 1)u=0 & \text{on } M,\\
				\partial_\nu u = \psi & \text{on } \partial M.
			\end{cases}
    \end{equation}
    The corresponding Neumann-to-Dirichlet map is defined by
    \begin{align*}
			& \mathcal{N}: H^{1/2}(\partial M)^* \mapsto H^{1/2}(\partial M),\\
			& \mathcal{N}\psi = u^{\psi}\arrowvert_{\partial M},
    \end{align*} 
    where $u^\psi$ is the solution to \eqref{hol_1}. Using Green's identity, we obtain 
    \begin{align*}
        \lambda_j(u_j,u^\psi)_{L^2(M,e^\phi d\mu_g)} &= - \int_M \Delta_g^Fu_j(x) u^\psi(x) e^{\phi(x)} d\mu_g(x)\\
        &= \int_M u_j(x) \Delta_g^Fu^\psi(x) e^{\phi(x)} d\mu_g(x) + \langle u_je^\phi\arrowvert_{\partial M}, \psi\rangle.
    \end{align*}
    Therefore,
    \begin{equation*}
        \langle u_je^\phi\arrowvert_{\partial M}, \psi\rangle = (\lambda_j - 1) (u_j,u^\psi)_{L^2(M,e^\phi d\mu_g)}.
    \end{equation*}
    Then,
    \begin{equation*}
        U(x) =  (\lambda_k - \lambda) \sum_{j\neq k} \frac{\lambda_j - 1}{(\lambda_j - \lambda_k)(\lambda_j - \lambda)}(u_j,u^\psi)_{L^2(M,e^\phi d\mu_g)} u_j(x) .
    \end{equation*}
    Let us set
    \begin{equation*}
        I_1(x) := \sum_{j\neq k} \frac{1}{\lambda_j - \lambda}(u_j,u^\psi)_{L^2(M,e^\phi d\mu_g)} u_j(x)
    \end{equation*}
    and
    \begin{equation*}
        I_2(x) := \sum_{j\neq k} \frac{\lambda_k - 1}{(\lambda_j - \lambda_k)(\lambda_j - \lambda)}(u_j,u^\psi)_{L^2(M,e^\phi d\mu_g)} u_j(x),
    \end{equation*}
    so that
    \begin{equation*}
        U(x) = (\lambda_k - \lambda) \left( I_1(x) + I_2(x)\right).
    \end{equation*}
    By the spectral theorem, we know that
    \begin{equation*}
        (-\Delta_N^F - \lambda)^{-1} u^\psi = \sum_{j=1}^\infty \frac{1}{\lambda_j - \lambda}(u_j,u^\psi)_{L^2(M,e^\phi d\mu_g)} u_j(x).
    \end{equation*}
    Hence,
    \begin{equation*}
        I_1(x) = (Id - P_k) (-\Delta_N^F - \lambda)^{-1} u^\psi,
    \end{equation*}
    where $Id$ is the identity operator and $P_k$ is the spectral projection to $\{u_k\}$. Therefore, 
    \begin{equation*}
        \|I_1\|_{H^1(M)} \lesssim \|\nabla_gu^\psi\|_{L^2(M,e^\phi d\mu_g)} + \|u^\psi\|_{L^2(M,e^\phi d\mu_g)}.
    \end{equation*}
    Furthermore, using the divergence theorem, we compute
    \begin{align}\label{ND_map_dua}
        0 &= \int_M (-\Delta_g^F u^\psi(x) + u^\psi(x)) u^\psi(x) e^{\phi(x)} d\mu_g(x) \\
        \nonumber &= \|\nabla_gu^\psi\|_{L^2(M,e^\phi d\mu_g)}^2 + \|u^\psi\|_{L^2(M,e^\phi d\mu_g)}^2 - \langle e^\phi \mathcal{N} \psi,\psi \rangle,
    \end{align}
    which implies that 
    \begin{equation*}
        \|I_1\|_{H^1(M)} \lesssim \sqrt{\langle e^\phi \mathcal{N} \psi,\psi \rangle}.
    \end{equation*}
    
    Next, we estimate $I_2$:
    \begin{multline*}
        \|I_2\|_{H^1(M)} \approx \|\nabla_gI_2\|_{L^2(M,e^\phi d\mu_g)} + \|I_2\|_{L^2(M,e^\phi d\mu_g)}\\
        \lesssim \left( \sum_{j\neq k} \left(\frac{(\lambda_k - 1)}{(\lambda_j - \lambda_k)(\lambda_j - \lambda)}\right)^2\lambda_j + \left(\frac{(\lambda_k - 1)}{(\lambda_j - \lambda_k)(\lambda_j - \lambda)}\right)^2 \right)^{\frac{1}{2}}(u_j,u^\psi)_{L^2(M,e^\phi d\mu_g)}.
    \end{multline*}
    Therefore,
    \begin{equation*}
        \|I_2\|_{H^1(M)} \lesssim \left(\sum_{j\neq k} \frac{1}{\lambda_j^3} \right)^{\frac{1}{2}}(u_j,u^\psi)_{L^2(M,e^\phi d\mu_g)} \lesssim \|u^\psi\|_{L^2(M,e^\phi d\mu_g)}
    \end{equation*}
    Here, as a consequence of Weyl's law \cite{Weyl}, we used that $\lambda_j \approx j^{\frac{2}{3}}$. See also \cite{BNR} for weighted Laplace operators. Due to \eqref{ND_map_dua}, we derive 
    \begin{equation*}
        \|I_2\|_{H^1(M)} \lesssim \sqrt{\langle e^\phi \mathcal{N} \psi,\psi \rangle},
    \end{equation*}
    and hence, 
    \begin{equation*}
        \|U\|_{H^1(M)} \lesssim |\lambda_k - \lambda|\sqrt{\langle e^\phi \mathcal{N} \psi,\psi \rangle}.
    \end{equation*}
    Since $-1 \notin \mathrm{spec} (-\Delta_N^F)$, it follows that $\mathcal{N}$ is a bounded operator from $H^{1/2}(\partial M)^*$ to $H^{1/2}(\partial M)$, see for instance \cite{BEHRNDT20155903}. Therefore, we conclude 
    \begin{equation*}
        \|R_{\lambda_k,\lambda}\psi\|_{H^{1/2}(M)} \leq \|U\|_{H^1(M)} \lesssim |\lambda_k - \lambda| \|\psi\|_{H^{1/2}(\partial M)^*}.
    \end{equation*}
    This completes the proof.

	\end{proof}
	
	\section{Proof of the main result}\label{proof_main_result}
	In this section, we prove our main result. We first begin with auxiliary lemmas which will be used subsequently. The following lemma is well-known in spectral theory. We state it here for the reader's convenience.
	\begin{lemma}[Glazman Lemma]\label{Glazman}
		Let $A$ be a lower-semibounded self-adjoint operator in a Hilbert space $\left(\mathcal{H}, \langle \cdot,\cdot\rangle\right)$ with corresponding closed sesquilinear form $a$ and form domain $\mathrm D(a)$. Then, it holds
		\begin{equation*}
			N(\lambda,A)=\sup \{\mathrm{dim} L\mid~ L \text{ subspace of } \mathrm D(a) \text{ s.th. } a(u,u)<\lambda\langle u,u\rangle \text{ for } u\in L \setminus \{ 0\} \},
		\end{equation*}
		where $N(\lambda,A)$ is the spectral distribution function.
	\end{lemma}
	
	The proofs of the following two lemmas are based on proofs of Lemma 3.1 and Theorem 3.2 in \cite{Rohleder}, respectively.
	\begin{lemma}\label{uniq_con}
		Let $\lambda \in \mathbb{R}$ and $u\in H^1(M)$ such that $-\Delta_g^F u = \lambda u$. Let $x_0\in \partial M$. If $u\left.\right|_{B_h(x_0,\varepsilon)} = 0$ and  $\partial_{\nu}u\left.\right|_{B_h(x_0,\varepsilon)} = 0$, then $u=0$ identically on $M$. (Recall that we consider $\varepsilon$ being smaller than injectivity radius.)
	\end{lemma}
	
	\begin{proof}
		Let us extend $M$ to a compact connected smooth Riemannian manifold $\widetilde{M} $ such that $\overline{\widetilde{M}\setminus M}$ is compact with non-empty interior and $\overline{\widetilde{M}\setminus M} \cap M = B_h(x_0,\varepsilon)$. Let $\widetilde{u}$ be the extension by zero of $u$ to $\widetilde{M}$. Let $\widetilde{F}$ be any smooth extension, up to the boundary, of $F$ to $\widetilde{M}$, so that we get a new weighted Laplacian $-\widetilde{\Delta}_{\widetilde{g}}^{\widetilde{F}}$. Since $u\left.\right|_{B_h(x_0,\varepsilon)} = 0$ and  $\partial_{\nu}u\left.\right|_{B_h(x_0,\varepsilon)} = 0$, it follows that $\widetilde{u} \in H^1(\widetilde{M})$, moreover $-\widetilde{\Delta}_{\widetilde{g}}^{\widetilde{F}} \widetilde{u} = \lambda \widetilde{u}$. Since $\widetilde{u}\left.\right|_{\widetilde{M}\setminus M} = 0$, the unique continuation arguments imply that $u=0$ identically on $M$.
	\end{proof}
	The previous lemma can be used to derive the following strict monotonicity principle, which will be used later.
	\begin{lemma}\label{str manaton}
		Assume that $0< \varepsilon_1 < \varepsilon_2$, then 
		\begin{equation*}
			\lambda_{j,\varepsilon_1} < \lambda_{j,\varepsilon_2}, \qquad j\in \mathbb{N}.
		\end{equation*}
		
\end{lemma}
\begin{proof}
	Since $\varepsilon_1 < \varepsilon_2$, it follows that $\mathrm{D}(a_{\varepsilon_2}) \subset \mathrm{D}(a_{\varepsilon_1})$, and hence, by Lemma \ref{Glazman}, we know that $\lambda_{j,\varepsilon_1} \leq \lambda_{j,\varepsilon_2}$. We are now required to show that the previous estimate is strict. 
	
	Let $\lambda = \lambda_{j,\varepsilon_2}$ and $\delta>0$ be sufficiently small such that
	\begin{equation*}
		(\lambda, \lambda+\delta) \cap \mathrm{spec}(-\Delta_{Mix,\varepsilon_i}^F) = \emptyset, \qquad \text{for } i=1,2.
	\end{equation*}
	We denote
	\begin{equation*}
		k = N\left(\lambda+\delta, -\Delta_{Mix,\varepsilon_2}^F\right) \quad \text{and} \quad L:= \textrm{span}\{u_{1,\varepsilon_2}, \cdots, u_{k,\varepsilon_2}\}.
	\end{equation*}
	Then $k\geq j$ and 
	\begin{equation}\label{ineq_eps_1}
		a_{\varepsilon_1}(u,u) = a_{\varepsilon_2}(u,u) < (\lambda +\delta) (u,u)_{L^2(M,e^\phi d\mu_g)}, \qquad \text{for } u\in L. 
	\end{equation}
	
	Since $\Gamma_{\varepsilon_2,a}\setminus\Gamma_{\varepsilon_1,a}$ has a non-empty interior in $\partial M$, by Lemma \ref{uniq_con}, it follows that if $v\in \mathrm{Ker}(-\Delta_{Mix,\varepsilon_1}^F - \lambda)$ then $v \notin L$. Therefore, we obtain
	\begin{equation}\label{dim_direct_sum}
		\mathrm{dim} \left(\mathrm{Ker}(-\Delta_{Mix,\varepsilon_1}^F - \lambda) \oplus L \right) = \mathrm{dim} \left(\mathrm{Ker}\left(-\Delta_{Mix,\varepsilon_1}^F - \lambda\right)\right) + \mathrm{dim} L.
	\end{equation}
	
	Furthermore, since $\mathrm{Ker}\left(-\Delta_{Mix,\varepsilon_1}^F - \lambda\right) \subset \mathrm{D}(-\Delta_{Mix,\varepsilon_1}^F)$ and $L\subset  \mathrm{D}(a_{\varepsilon_2}) \subset \mathrm{D}(a_{\varepsilon_1}) $, Theorem 2.1 in \cite[Chapter 6]{Kato} and estimate \eqref{ineq_eps_1} give
	\begin{align*}
		a_{\varepsilon_1}(v + u, v + u) &< (\lambda +\delta) (v,v)_{L^2(M,e^\phi d\mu_g)} + (\lambda +\delta) (u,u)_{L^2(M,e^\phi d\mu_g)} + 2 a_{\varepsilon_1}(v , u)\\
		&< (\lambda +\delta) \left((v,v)_{L^2(M,e^\phi d\mu_g)} + (u,u)_{L^2(M,e^\phi d\mu_g)} \right) + 2\lambda (v,u)_{L^2(M,e^\phi d\mu_g)}\\
		& \leq (\lambda +\delta) (v + u, v + u)_{L^2(M,e^\phi d\mu_g)},
	\end{align*}
	for $v \in \mathrm{Ker}\left(-\Delta_{Mix,\varepsilon_1}^F - \lambda\right)$ and $u\in L$. Therefore, by Lemma \ref{Glazman}, it follows 
	\begin{equation*}
		N\left(\lambda+\delta, -\Delta_{Mix,\varepsilon_1}^F\right) \geq \mathrm{dim} \left(\mathrm{Ker}\left(-\Delta_{Mix,\varepsilon_1}^F - \lambda\right) \oplus L \right), 
	\end{equation*}
	and hence, by \eqref{dim_direct_sum},
	\begin{equation*}
		N\left(\lambda+\delta, -\Delta_{Mix,\varepsilon_1}^F\right) \geq \mathrm{dim} \left(\mathrm{Ker}\left(-\Delta_{Mix,\varepsilon_1}^F - \lambda\right)\right) + \mathrm{dim} L.
	\end{equation*}
	Then
	\begin{equation*}
		N\left(\lambda, -\Delta_{Mix,\varepsilon_1}^F\right) = N\left(\lambda+\delta, -\Delta_{Mix,\varepsilon_1}^F\right) - \mathrm{dim} \left(\mathrm{Ker}\left(-\Delta_{Mix,\varepsilon_1}^F - \lambda\right)\right) \geq k \geq j,
	\end{equation*}
	so that $\lambda_{j,\varepsilon_1} < \lambda_{j,\varepsilon_2}$.
\end{proof}

Since $\mathrm{D}(a_{\varepsilon}) \subset \mathrm{D}(a_{N})$, it follows that  $\lambda_{j,\varepsilon}$ is bounded from below by $\lambda_j$ and decreases as $\varepsilon\rightarrow 0$, by the last lemma. Therefore, we can define the following limit 
\begin{equation*}
	\lambda_{j,0} := \lim_{\varepsilon\rightarrow 0} \lambda_{j,\varepsilon}.
\end{equation*}
Next, we show that $\{\lambda_{j,0}\}_{j\in\mathbb{N}}$ coincides with the sequence of eigenvalues of $-\Delta_N^F$: 
\begin{lemma}\label{biject}
	For any $j\in\mathbb{N}$, the equality $\lambda_{j,0} = \lambda_j$ holds. 
\end{lemma}
\begin{proof}
	For $\varepsilon>0$, we know that $\mathrm{D}(a^N) \subset \mathrm{D}(a_\varepsilon)$. Therefore, by Lemma \ref{Glazman}, $\lambda_j\leq \lambda_{j,\varepsilon}$. Recalling our definition of $\lambda_{j,0}$, we conclude that $\lambda_j\leq \lambda_{j,0}$, or equivalently
	\begin{equation*}
		N(\lambda, -\Delta_N^F) \geq \#\{\lambda_{j,0}: \; \lambda_{j,0}<\lambda\}, \qquad \lambda>0.
	\end{equation*}
	Therefore, to prove $\lambda_{j,0} = \lambda_j$, it suffices to show that if $\lambda\in \mathrm{spec} (-\Delta_N^F)$ with multiplicity $l$, then $\lambda$ appears in $\{\lambda_{j,0}\}_{j\in \mathbb{R}}$ at least $l$ times.
	
	Let $\lambda\in \mathrm{spec} (-\Delta_N^F)$ and $l$ be its multiplicity. Then there exists $k\in \mathbb{N}$ such that $\lambda < \lambda_{k+1}$ and 
	\begin{equation*}
		\lambda_{k - l +1}=\cdots = \lambda_k = \lambda.
	\end{equation*}
	Therefore, there exists $\alpha_0>0$ such that $N(\lambda + \alpha, -\Delta_N^F) = k$ for any $\alpha\in (0,\alpha_0)$. For any $\alpha\in (0,\alpha_0)$, we aim to find a small $\varepsilon>0$ so that $N(\lambda + \alpha, -\Delta_{Mix,\varepsilon}^F) = k$. 
	
	Let $\chi_\varepsilon\in C^\infty(M)$ denote a smooth cutoff function such that 
	\begin{equation*}
		\chi_\varepsilon(x)=
		\begin{cases}
			1 & \text{for } x\in M\setminus B_g(x^*,3\varepsilon),\\
			0 & \text{for } x\in B_g(x^*,2\varepsilon)
		\end{cases}
	\end{equation*}
	and 
	\begin{equation*}
		\|\nabla_g\chi_\varepsilon\|_{L^2(M,e^\phi d\mu_g)} \rightarrow 0, \text{for } \varepsilon\rightarrow 0.
	\end{equation*}
	Consider the set of functions
	\begin{equation*}
		L_{\varepsilon}:= \{u_{1}\chi_{\varepsilon}, \cdots, u_{k}\chi_{\varepsilon}\}.
	\end{equation*}
	Since $\{u_j\}_{j=1}^k$ are linearly independent in $L^2$ and $u_{j}\chi_{\varepsilon} \rightarrow u_{j}$ in $L^2$, it follows that $\{u_j\chi_{\varepsilon}\}_{j=1}^k$ are also linearly independent in $L^2$ for sufficiently small $\varepsilon>0$, so that $\mathrm{dim}(L_\varepsilon) = k$. The definition of $\chi_{\varepsilon}$ implies that $L_{\varepsilon} \subset \mathrm{D}(a_\varepsilon)$ and 
	\begin{equation*}
		a(u_j\chi_{\varepsilon},u_j\chi_{\varepsilon}) \rightarrow a(u_j,u_j),
	\end{equation*}
	\begin{equation*}
		(u_j\chi_{\varepsilon},u_j\chi_{\varepsilon})_{L^2(M,e^\phi d\mu_g)} \rightarrow (u_j,u_j)_{L^2(M,e^\phi d\mu_g)},
	\end{equation*}
	as $\varepsilon\rightarrow 0$. Therefore, for $\alpha \in (0,\alpha_0)$, there exists $\varepsilon>0$ such that 
	\begin{equation*}
		\frac{a(u,u)}{(u,u)_{L^2(M,e^\phi d\mu_g)}} < \lambda + \alpha \qquad \text{for } u\in L_{\varepsilon},
	\end{equation*}
	which implies that $N(\lambda +\alpha, -\Delta_{Mix,\varepsilon}^F) = k$. Therefore, $\mathrm{spec}(-\Delta_{Mix,\varepsilon}^F)\cap (\lambda,\lambda+\alpha)$ is not empty. Moreover, since $\lambda_j<\lambda_{j,\varepsilon}$, we conclude 
	\begin{equation*}
		\left\{ \lambda_{k - l +1,\varepsilon},\cdots, \lambda_{k,\varepsilon}\right\} \subset \mathrm{spec}(-\Delta_{Mix,\varepsilon}^F)\cap (\lambda,\lambda+\alpha).
	\end{equation*}
	Since $\alpha_0$ is an arbitrary sufficiently small number, we conclude that $\lambda$ appears in $\{\lambda_{j,0}\}_{j\in \mathbb{N}}$ at least $l$ times. This completes the proof.
\end{proof}

Next, we show that the eigenfunctions of $-\Delta_{Mix,\varepsilon}^F$ and $-\Delta_N^F$ are close to each other in the following sense:
\begin{lemma}\label{ef}
	Assume that $\lambda_{j}$ is a simple eigenvalue of $-\Delta_N^F$. Then there exists $C>0$ such that
	\begin{equation*}
		\left|(u_j,u_{j, \varepsilon})_{L^2(M,e^\phi d\mu_g)}\right|>C,
	\end{equation*}
	for sufficiently small $\varepsilon>0$.
\end{lemma}
\begin{proof}
	Recall that $\{u_k\}_{k\in \mathbb{N}}$ forms an orthonormal basis on $L^2(M,e^\phi d\mu_g)$, so that we can express
	\begin{equation*}
		u_{j,\varepsilon}(x) = \sum_{k\in\mathbb{N}} c_k^\varepsilon u_k.
	\end{equation*}
	Assume that the lemma is false, which means that there is a sequence of positive numbers $\{\varepsilon_l\}_{l\in \mathbb{N}}$ such that 
	\begin{equation}\label{assumpepsilon}
		\varepsilon_l \rightarrow 0, \qquad c_j^{\varepsilon_l} \rightarrow 0,
	\end{equation}
	as $l\rightarrow\infty$. Since $\lambda_j$ is simple we can choose $\alpha>0$ such that $\lambda_j + \alpha < \lambda_{j+1}$.  Let us define
	\begin{equation*}
		\omega_{j,\varepsilon_l} := \sum_{k\neq j} c_k^{\varepsilon_l} u_k.
	\end{equation*}
    Then
    \begin{multline*}
        (\omega_{j,\varepsilon_l} , \omega_{j,\varepsilon_l} )_{L^2(M,e^\phi d\mu_g)} = (u_{j,\varepsilon_l} - c_j^{\varepsilon_l}u_j, u_{j,\varepsilon_l} - c_j^{\varepsilon_l}u_j)_{L^2(M,e^\phi d\mu_g)}\\
        =(u_{j,\varepsilon_l}, u_{j,\varepsilon_l})_{L^2(M,e^\phi d\mu_g)} - 2(u_{j,\varepsilon_l} , c_j^{\varepsilon_l}u_j)_{L^2(M,e^\phi d\mu_g)} + (c_j^{\varepsilon_l}u_j, c_j^{\varepsilon_l}u_j)_{L^2(M,e^\phi d\mu_g)},
    \end{multline*}
    and hence,
    \begin{equation*}
        (\omega_{j,\varepsilon_l} , \omega_{j,\varepsilon_l} )_{L^2(M,e^\phi d\mu_g)} = (u_{j,\varepsilon_l}, u_{j,\varepsilon_l})_{L^2(M,e^\phi d\mu_g)} + o(1)
    \end{equation*}
    as $l\rightarrow\infty$. Similarly, we obtain
    \begin{multline}\label{quad_form_omega}
        a^N(\omega_{j,\varepsilon_l} , \omega_{j,\varepsilon_l} ) = a^N (u_{j,\varepsilon_l}, u_{j,\varepsilon_l}) - 2 a^N (u_{j,\varepsilon_l}, c_j^{\varepsilon_l}u_j) + a^N (c_j^{\varepsilon_l}u_j, c_j^{\varepsilon_l}u_j)\\
        = \lambda_{j\varepsilon_l} (u_{j,\varepsilon_l}, u_{j,\varepsilon_l})_{L^2(M,e^\phi d\mu_g)} - \lambda_j(c_j^{\varepsilon_l})^2 = \lambda_{j\varepsilon_l} (\omega_{j,\varepsilon_l}, \omega_{j,\varepsilon_l})_{L^2(M,e^\phi d\mu_g)} + o(1)
    \end{multline}
    as $l\rightarrow\infty$. Since $\lambda_{j,\varepsilon_l} \rightarrow \lambda_j$ as $l\rightarrow \infty$, this implies that
    \begin{equation*}
		a^N(\omega_{j,\varepsilon_l} , \omega_{j,\varepsilon_l} )  < (\lambda_j + \alpha)  (\omega_{j,\varepsilon_l} , \omega_{j,\varepsilon_l} )_{L^2(M,e^\phi d\mu_g)}
    \end{equation*}
    for sufficiently large $l\in \mathbb{N}$. Let us show that 
    \begin{equation}\label{omega_not_in_span}
		\omega_{j,\varepsilon_l} \notin \text{span}\{u_1,\cdots,u_j\}
    \end{equation}
    for sufficiently large $l\in \mathbb{N}$. Assume this is not true. Without loss of generality, we assume that \eqref{omega_not_in_span} is false for all $l\in \mathbb N$, otherwise consider a subsequence. Due to the definition of $\omega_{\omega,\varepsilon_l}$, this implies that
    \begin{equation*}
		\omega_{j,\varepsilon_l} \notin \text{span}\{u_1,\cdots,u_{j-1}\}.
    \end{equation*}
    In this case, we would have
    \begin{equation*}
        a^N(\omega_{j,\varepsilon_l} , \omega_{j,\varepsilon_l} ) \leq \lambda_{j-1}(\omega_{j,\varepsilon_l} , \omega_{j,\varepsilon_l} )_{L^2(M,e^\phi d\mu_g)}.
    \end{equation*}
    Since $\lambda_j$ is simple, this contradicts to \eqref{quad_form_omega}. Therefore, \eqref{omega_not_in_span} holds.
    
    We now let $u\in \text{span}\{u_1,\cdots,u_j\}$. Then 
	\begin{multline}\label{omega_plus_u}
		a^N(\omega_{j,\varepsilon_l} + u , \omega_{j,\varepsilon_l} + u)\\
             \leq(\lambda_{j} + \alpha) (\omega_{j,\varepsilon_l} ,\omega_{j,\varepsilon_l} )_{L^2(M,e^\phi d\mu_g)} + \lambda_{j} (u ,u )_{L^2(M,e^\phi d\mu_g)} + 2a^N(\omega_{j,\varepsilon_l} ,u),
	\end{multline}
    for sufficiently large $l\in \mathbb N$. Let us estimate, the last term of the right-hand side. Let $\chi_{\varepsilon_l}$ be the function described in the proof of Lemma \ref{biject}, then
    \begin{multline*}
        a^N(\omega_{j,\varepsilon_l} ,u) = a^N(u_{j,\varepsilon_l},u) - c_j^{\varepsilon_l}a^N(u_j,u)\\
        = a^N(u_{j,\varepsilon_l},\chi_{\varepsilon_l}u) + a^N(u_{j,\varepsilon_l},u -  \chi_{\varepsilon_l}u) - c_j^{\varepsilon_l}a^N(u_j,u) = a^N(u_{j,\varepsilon_l},\chi_{\varepsilon_l}u) +o(1)
    \end{multline*}
    as $l\rightarrow \infty$. Since $u_{j,\varepsilon_l} \in \mathrm{D}(-\Delta_{Mix,\varepsilon_l}^N)$ and $u\chi_{\varepsilon} \in \mathrm{D}(a_{\varepsilon_l})$, it follows that
    \begin{equation*}
        a^N(u_{j,\varepsilon_l},\chi_{\varepsilon_l}u) = \lambda_{j,\varepsilon_l}(u_{j,\varepsilon_l},\chi_{\varepsilon_l}u)_{L^2(M,e^\phi d\mu_g)} = \lambda_{j,\varepsilon_l}(\omega_{j,\varepsilon_l},u)_{L^2(M,e^\phi d\mu_g)} + o(1),
    \end{equation*}
    as $l\rightarrow \infty$. Therefore,
    \begin{equation*}
        a^N(\omega_{j,\varepsilon_l} ,u) \leq (\lambda_j + \alpha) (\omega_{j,\varepsilon_l},u)_{L^2(M,e^\phi d\mu_g)}.
    \end{equation*}
    
    Therefore, \eqref{omega_plus_u} gives
	\begin{equation*}
		a^N(\omega_{j,\varepsilon_l} + u , \omega_{j,\varepsilon_l} + u) \leq (\lambda_j + \alpha) (\omega_{j,\varepsilon_l} + u , \omega_{j,\varepsilon_l} + u)_{L^2(M,e^\phi d\mu_g)}, 
	\end{equation*}
    for sufficiently large $l\in \mathbb{N}$. Due to \eqref{omega_not_in_span} and Lemma \ref{Glazman}, we obtain
    $$N(\lambda_j + \alpha, \Delta_N^F) \geq j + 1.$$
    This contradicts to $\lambda_j + \alpha < \lambda_{j + 1}$.
\end{proof}

Now, we are ready to prove the main result.

\begin{proof}[Proof of Theorem \ref{main_for_ellipse}]
	Let $V_j\subset \mathbb{C}$ be an open neighbourhood of $\lambda_j$ which does not contain any other eigenvalues of $-\Delta_N^F$. Since $\lambda_j$ is simple, Theorem \ref{biject} implies that, for sufficiently small $\varepsilon>0$, $\lambda_{j,\varepsilon}$ is the only eigenvalue of $-\Delta_{Mix,\varepsilon}^F$ in $V_j$. For $\omega\in V_j$ and $x$, $y\in \partial M$, Proposition \ref{sing_stract} gives
	\begin{align*}
		G^{\omega}_{\partial M}(x,y)=&\frac{1}{2\pi}d_g(x,y)^{-1} -\frac{H(x)}{4\pi}\log d_h(x,y)+\frac{g_x(F,\nu)}{4\pi}\log d_h(x,y)\\ &+\frac{1}{16\pi}\left(\mathrm{II}_{x}\left(\frac{\exp_{x}^{-1}(y)}{|\exp_{x}^{-1}(y)|_h}\right)-\mathrm{II}_{x}\left(\frac{\star\exp_{x}^{-1}(y)}{|\exp_{x}^{-1}(y)|_h}\right)\right) \\
		& +\frac{1}{4\pi}h_{x}\left(F^{||}(x),\frac{\exp^{-1}_{x}(y)}{|\exp^{-1}_{x}(y)|_h}\right)\\
		& + \frac{u_j(x)u_j(y)}{\lambda_j - \omega^2}e^{\phi(y)}  + R_{\partial M}^{\lambda_{j,\varepsilon}}(x,y).
	\end{align*}
	From Green's identity, we know that
	\begin{equation*}
		u_{j,\varepsilon} = (\lambda_{j,\varepsilon} - \omega^2) \int_{M} G_{M}^\omega(x,y) u_{j,\varepsilon}(y) d\mu_g(y) + \int_{\Gamma_{\varepsilon, a}} G_M^{\omega}(x,y) \partial_{\nu} u_{j,\varepsilon}(y) d\mu_h(y).
	\end{equation*}
	We choose $\omega^2 = \lambda_{j,\varepsilon}$ and restrict the last identity to $\Gamma_{\varepsilon,a}$, to obtain
	\begin{equation}\label{gpu0}
		\int_{\Gamma_{\varepsilon, a}} G_{\partial M}^{\sqrt{\lambda_{j,\varepsilon}}}(x,y) \partial_{\nu} u_{j,\varepsilon}(y)d\mu_h(y) = 0.
	\end{equation}
	Next, we will use the coordinate system given by \eqref{res coord}. We denote 
	\begin{equation*}
		\tilde{u}_j(t'): = u_j(x^{\varepsilon}(t_1,at_2)),  \quad \tilde{\phi}(t'): = \phi(x^{\varepsilon}(t_1,at_2))
	\end{equation*}
	and
	\begin{equation*}
		v_\varepsilon : = \partial_{\nu} u_{j,\varepsilon},\quad \tilde{v}_\varepsilon(t'): = v_\varepsilon(x^{\varepsilon}(t_1,at_2)).
	\end{equation*}
	Note that in these coordinates the volume form for $\partial M$ is given by
	\begin{equation}\label{volform}
		d\mu_h(y) = a \varepsilon^2 (1 + \varepsilon^2Q_\varepsilon(s')) ds_1 \wedge ds_2,
	\end{equation}
	for some smooth function $Q_\varepsilon$ whose derivatives of all orders are bounded uniformly in $\varepsilon$. Therefore, if we put the expression for $G^{\sqrt{\lambda_{j,\varepsilon}}}_{\partial M}$ into \eqref{gpu0} and use Lemmas \ref{l1}-\ref{l4}, we obtain
	\begin{align}\label{after coord system}
		0 =&\frac{1}{2\pi} \varepsilon L_a\tilde{v}_\varepsilon + a\varepsilon^2\frac{\tilde{u}_j(t')}{\lambda_j - \lambda_{j,\varepsilon}} \int_{\mathbb{D}} \tilde{u}_j(s')\tilde{v}_\varepsilon(s')e^{\tilde{\phi}(s')}ds'\\
		\nonumber &  - \frac{1}{4\pi} \varepsilon^2 (H(x^*)-\partial_{\nu}\phi(x^*))R_{\log,a} \tilde{v}_\varepsilon + \frac{1}{16\pi} \varepsilon^2(\kappa_1(x^*) - \kappa_2(x^*)) R_{\infty,a} \tilde{v}_\varepsilon \\
		\nonumber 
        & + \frac{1}{4\pi} \varepsilon ^2 R_{F,a}\tilde{v}_\varepsilon + \varepsilon^2 R_a^{\lambda_{j,\varepsilon}}\tilde{v}_\varepsilon - \frac{1}{4\pi} \varepsilon^2\log\varepsilon(H(x^*)-\partial_{\nu}\phi(x^*)) R_{I,a}\tilde{v}_\varepsilon + \varepsilon^3\log\varepsilon\mathcal{A}_\varepsilon\tilde{v}_{\varepsilon},
	\end{align}
	where $R_a^{\omega}: C_c^\infty(\mathbb{D}) \mapsto D'(\mathbb{D})$ is the operator given by the kernel $$aR_{\partial M}^{\omega}(x^{\varepsilon}(t_1,at_2),x^{\varepsilon}(s_1,as_2))$$ 
	for $\omega \in V_j$, and $\A_\varepsilon : H^{1/2}(\D; ds')^* \mapsto H^{1/2}(\D; ds')$ is an operator with the norm bounded uniformly in $\varepsilon$. From now on, we will denote by $\A_\varepsilon$ any operator which takes $H^{1/2}(\D ;ds')^* \mapsto H^{1/2}(\D; ds')$ whose operator norm is bounded uniformly in $\varepsilon$.
	
	Let us denote 
	\begin{equation*}
		\mathcal{R}^{\lambda_j,\lambda_{j,\varepsilon}} := \varepsilon^2 R_a^{\lambda_{j,\varepsilon}} - \varepsilon^2 R_a^{\lambda_{j}} 
	\end{equation*}
	\begin{multline*}
		\mathcal{R}_\varepsilon := - \frac{1}{4\pi} \varepsilon^2\log\varepsilon(H(x^*)-\partial_{\nu}\phi(x^*)) R_{I,a} - \frac{1}{4\pi} \varepsilon^2 (H(x^*)-\partial_{\nu}\phi(x^*))R_{\log,a} \\
		+ \frac{1}{16\pi} \varepsilon^2(\kappa_1(x^*) - \kappa_2(x^*)) R_{\infty,a} + \frac{1}{4\pi} \varepsilon ^2 R_{F,a} + \varepsilon^2 R^{\lambda_j}_{\partial M}(x^*,x^*)R_{I,a} + \varepsilon^3\log\varepsilon\mathcal{A}_\varepsilon.
	\end{multline*}
	By Lemma 5.1 in \cite{NursultanovTzouTzou}, we know that 
	\begin{equation*}
		\left\| R_a^{\lambda_{j}}  -  R^{\lambda_j}_{\partial M}(x^*,x^*)R_{I,a}\right\|_{H^{1/2}(\D; ds')^* \mapsto H^{1/2}(\D; ds')} = O(\varepsilon \log \varepsilon),
	\end{equation*}
	and hence, \eqref{after coord system} becomes
	\begin{equation}\label{LUR}
		0 = \frac{1}{2\pi} \varepsilon L_a\tilde{v}_\varepsilon + a\varepsilon^2\frac{\tilde{u}_j(t')}{\lambda_j - \lambda_{j,\varepsilon}} \int_{\mathbb{D}} \tilde{u}_j(s')\tilde{v}_\varepsilon(s')e^{\tilde{\phi}(s')}ds' + \mathcal{R}_\varepsilon \tilde{v}_\varepsilon + \mathcal{R}^{\lambda_j,\lambda_{j,\varepsilon}}\tilde{v}.
	\end{equation}
	Assume that
	\begin{equation}\label{assumption psi zero}
		\int_{\mathbb{D}} \tilde{u}_j(s')\tilde{v}_\varepsilon(s')e^{\tilde{\phi}(s')}ds' = 0.
	\end{equation}
	Then, recalling \eqref{volform}, we get
	\begin{equation*}
		\int_{\Gamma_{\varepsilon, a}} u_j(x) \partial_{\nu}u_{j,\varepsilon}(x) e^{\phi(x)}d\mu_h(x) = O(\varepsilon^4).
	\end{equation*}
	Using Green's identity, we derive
	\begin{multline*}
		(\lambda_{j,\varepsilon} - \lambda_j)(u_j, u_{j,\varepsilon})_{L^2(M,e^\phi d\mu_g)} = \int_{M} (\Delta_g^F u_j(x)u_{j,\varepsilon}(x) -  u_j(x) \Delta_g^F u_{j,\varepsilon}(x)) e^{\phi(x)} d\mu_g(x) \\
		= - \int_{\Gamma_{\varepsilon, a}} u_j(x) \partial_\nu u_{j,\varepsilon}(x) e^{\phi(x)} d\mu_g(x) = O(\varepsilon^4).
	\end{multline*}
	Then, by Lemma \ref{ef}, it follows 
	\begin{equation*}
		|\lambda_{j,\varepsilon} -  \lambda_j| = O(\varepsilon^4).
	\end{equation*}
	Therefore, by using Proposition \ref{tr} and recalling Remark \ref{rem}, we estimate
	\begin{equation*}
		\|\mathcal{R}_\varepsilon + \mathcal{R}^{\lambda_j,\lambda_{j,\varepsilon}}\|_{H^{1/2}(\mathbb{D})^* \mapsto H^{1/2}(\mathbb{D})} = O(\varepsilon^2\log\varepsilon).
	\end{equation*}
	Since  $L_a$ is invertable as an operator from $H^{1/2}(\mathbb{D})^*$ to $H^{1/2}(\mathbb{D})$, see Section 4 in \cite{NursultanovTzouTzou}, it follows that
	\begin{equation*}
		L_a + \frac{2\pi}{\varepsilon}\left( \mathcal{R}_\varepsilon + \mathcal{R}^{\lambda_j,\lambda_{j,\varepsilon}}\right): H^{1/2}(\mathbb{D})^* \mapsto H^{1/2}(\mathbb{D})
	\end{equation*}
	is an invertable operator. Therefore, \eqref{LUR} and \eqref{assumption psi zero} imply that $\tilde{v}_\varepsilon = 0$ on $\mathbb{D}$, and hence, $\partial_{\nu} u_{j,\varepsilon} = 0$ on $\Gamma_{\varepsilon, a}$. Then, by Lemma \ref{uniq_con}, we would have $u_{j,\varepsilon} = $ on $M$, and hence $\lambda_{j,0} = 0$. This contradicts to $\lambda_{j,\varepsilon}>\lambda_j\geq 0$. Therefore,
	\begin{equation*}
		\int_{\mathbb{D}} \tilde{u}_j(s')\tilde{v}_\varepsilon(s')e^{\tilde{\phi}(s')}ds' \neq 0.
	\end{equation*}
	Therefore, we can define
	\begin{equation*}
		\tilde{\psi}_\varepsilon : = \frac{\tilde{v}_\varepsilon}{\int_{\mathbb{D}} \tilde{u}_j(s')\tilde{v}_\varepsilon(s')e^{\tilde{\phi}(s')}ds'}.
	\end{equation*}
	Then, \eqref{LUR} becomes
	\begin{equation*}
		0 = \frac{1}{2\pi} \varepsilon L_a\tilde{\psi}_\varepsilon + a\varepsilon^2\frac{\tilde{u}_j(t')}{\lambda_j - \lambda_{j,\varepsilon}}  + \left(\mathcal{R}_\varepsilon + \mathcal{R}^{\lambda_j,\lambda_{j,\varepsilon}} \right)\tilde{\psi}_\varepsilon.
	\end{equation*}
	Let us hit both sides by $\frac{2\pi}{\varepsilon} L_a^{-1}$, to obtain
	\begin{equation}\label{exp_for_psi_tilde}
		0 = \tilde{\psi}_{\varepsilon} + 2\pi a\varepsilon\frac{L_{a}^{-1}\tilde{u}_j}{\lambda_j - \lambda_{j,\varepsilon}}  + \frac{2\pi}{\varepsilon}L_{a}^{-1}\left( \mathcal{R}_\varepsilon + \mathcal{R}^{\lambda_j,\lambda_{j,\varepsilon}}\right) \tilde{\psi}_\varepsilon.
	\end{equation}
    Since 
    \begin{equation*}
        \left\| L_a^{-1}\left( \mathcal{R}_\varepsilon + \mathcal{R}^{\lambda_j,\lambda_{j,\varepsilon}}\right) \right\|_{H^{1/2}(\mathbb{D})^* \mapsto H^{1/2}(\mathbb{D})} = O(\varepsilon^2\log\varepsilon),
    \end{equation*}
    relation \eqref{exp_for_psi_tilde}, implies that 
    \begin{equation}\label{est_for_psi}
        \|\tilde{\psi}\|_{H^{1/2}(\mathbb D)^*} = \frac{1}{\lambda_j - \lambda_{j,\varepsilon}}O(\varepsilon).
    \end{equation}
    This will be used later. Now, we multiply \eqref{exp_for_psi_tilde} by $\tilde{u}_je^{\tilde{\phi}}$ and integrate over $\mathbb{D}$ to derive
	\begin{equation*}
		0 = 1 + 2\pi a\varepsilon\frac{\langle L_{a}^{-1}[\tilde{u}_j], \tilde{u}_je^{\tilde{\phi}}\rangle}{\lambda_j - \lambda_{j,\varepsilon}}  + \frac{2\pi}{\varepsilon}\langle L_{a}^{-1}\left( \mathcal{R}_\varepsilon + \mathcal{R}^{\lambda_j,\lambda_{j,\varepsilon}}\right) \tilde{\psi}_\varepsilon, \tilde{u}_j e^{\tilde{\phi}}\rangle.
	\end{equation*}
	Equivalently, we write 
	\begin{equation*}
		\lambda_{j,\varepsilon} - \lambda_j = 2\pi a\varepsilon \langle L_{a}^{-1}\tilde{u}_j, \tilde{u}_je^{\tilde{\phi}}\rangle + \frac{2\pi}{\varepsilon}(\lambda_j - \lambda_{j,\varepsilon})\langle L_{a}^{-1}\left( \mathcal{R}_\varepsilon + \mathcal{R}^{\lambda_j,\lambda_{j,\varepsilon}}\right) \tilde{\psi}_\varepsilon, \tilde{u}_je^{\tilde{\phi}}\rangle.
	\end{equation*}
	Let us put \eqref{exp_for_psi_tilde} into the equation above to obtain
	\begin{multline}\label{lminl}
		\lambda_{j,\varepsilon} - \lambda_j = 2\pi a\varepsilon \langle L_{a}^{-1}\tilde{u}_j, \tilde{u}_je^{\tilde{\phi}}\rangle - 4\pi^2 a \langle L_{a}^{-1}\left( \mathcal{R}_\varepsilon + \mathcal{R}^{\lambda_j,\lambda_{j,\varepsilon}}\right) L_{a}^{-1}\tilde{u}_j, \tilde{u}_je^{\tilde{\phi}}\rangle\\
		- \frac{4\pi^2}{\varepsilon^2} (\lambda_j - \lambda_{j,\varepsilon})\langle L_{a}^{-1}\left( \mathcal{R}_\varepsilon + \mathcal{R}^{\lambda_j,\lambda_{j,\varepsilon}}\right) L_{a}^{-1}\left( \mathcal{R}_\varepsilon + \mathcal{R}^{\lambda_j,\lambda_{j,\varepsilon}}\right) \tilde{\psi}_\varepsilon, \tilde{u}_je^{\tilde{\phi}}\rangle.
	\end{multline}
    Taking into account \eqref{est_for_psi}, the last identity implies that 
    \begin{equation*}
        \lambda_{j,\varepsilon} - \lambda_\varepsilon = O(\varepsilon),
    \end{equation*}
    so that Proposition \ref{tr} gives
	\begin{equation*}
		\left\| \mathcal{R}^{\lambda_j,\lambda_{j,\varepsilon}}\right\|_{H^{1/2}(\mathbb{D})^* \mapsto H^{1/2}(\mathbb{D})} = O(\varepsilon^4).
	\end{equation*}
    Hence, we can put $\mathcal{R}^{\lambda_j,\lambda_{j,\varepsilon}}$ into $\varepsilon^2\log\varepsilon\mathcal{A}_{\varepsilon}$ term in the definition of $\mathcal{R}_\varepsilon$. Further, from \eqref{est_for_psi}, we obtain
	\begin{equation*}
		\langle L_{a}^{-1}\mathcal{R}_\varepsilon L_{a}^{-1}\mathcal{R}_\varepsilon\tilde{\psi}_\varepsilon, \tilde{u}_je^{\tilde{\phi}}\rangle = \frac{1}{\lambda_j - \lambda_{j,\varepsilon}}O\left(\varepsilon^5\log^2\varepsilon\right).
	\end{equation*}
	Therefore, equation \eqref{lminl} gives
	\begin{equation*}
		\lambda_{j,\varepsilon} = \lambda_j + 2\pi\varepsilon a \langle L_{a}^{-1}\tilde{u}_j, \tilde{u}_je^{\tilde{\phi}} \rangle - 4\pi^2 a \langle L_{a}^{-1}\mathcal{R}_\varepsilon L_{a}^{-1}\tilde{u}_j, \tilde{u}_je^{\tilde{\phi}} \rangle + O(\varepsilon^3\log^2\varepsilon).
	\end{equation*}
	Recalling the definition of $\mathcal{R}_{\varepsilon}$ and the boundedness of $\A_\varepsilon : H^{1/2}(\D; ds')^* \to H^{1/2}(\D; ds')$ we obtain
	\begin{align}\label{lambda minus lambda 1}
		\nonumber\lambda_{j,\varepsilon} - \lambda_j =& 2\pi\varepsilon a \int_{\mathbb{D}} L_a^{-1}\tilde{u}_j(s') \tilde{u}_j(s') e^{\tilde{\phi}(s')} ds'\\
		\nonumber & + \varepsilon^2\log\varepsilon a\pi (H(x^*)-\partial_{\nu}\phi(x^*)) \langle L_a^{-1} R_{I,a} L_a^{-1}\tilde{u}_j, \tilde{u}_je^{\tilde{\phi}} \rangle\\
		&+ \varepsilon^2 a \pi  (H(x^*)-\partial_{\nu}\phi(x^*)) \langle L_a^{-1} R_{log,a} L_a^{-1}\tilde{u}_j, \tilde{u}_je^{\tilde{\phi}} \rangle\\
		\nonumber& - \varepsilon^2 a \frac{\pi}{4} \left(\kappa_1(x^*) -  \kappa_2(x^*)\right) \langle L_a^{-1} R_{\infty,a} L_a^{-1}\tilde{u}_j, \tilde{u}_je^{\tilde{\phi}} \rangle\\
		\nonumber& - \varepsilon^2 a \pi \langle L_a^{-1} R_{F,a} L_a^{-1}\tilde{u}_j, \tilde{u}_je^{\tilde{\phi}} \rangle\\
		\nonumber& - \varepsilon^2  4\pi^2 a R_{\partial M}^{\lambda_j}(x^*,x^*)\langle L_a^{-1}R_{I,a} L_a^{-1}\tilde{u}_j, \tilde{u}_je^{\tilde{\phi}}\rangle\\
		\nonumber& + O(\varepsilon^3\log^2\varepsilon). 
	\end{align}
	We recall the definitions of $\tilde u_j$, $\tilde \phi$ and use Taylor series, to obtain
	\begin{multline*}
		\int_{\mathbb{D}} L_a^{-1}\tilde{u}_j(s') \tilde{u}_j(s') e^{\tilde{\phi}(s')} ds' = \int_{\mathbb{D}} L_a^{-1}\tilde{u}_j(s') u_j(x(\varepsilon s_1,a\varepsilon s_2)) e^{\phi(x(\varepsilon s_1,a\varepsilon s_2))} ds'\\
		= \int_{\mathbb{D}} L_a^{-1}\tilde{u}_j(s') \left( u_j(x^*)e^{\phi(x^*)} + \varepsilon(c_1s_1 + c_2s_2) + R_2^1(\varepsilon s') \right) ds'
	\end{multline*}
	where $R_2^1$ is the reminder term of the Taylor series for $\tilde u_j e^{\tilde \phi}$ near zero and $c_1$, $c_2$ are appropriate constants.
	Using (4.4) in \cite{NursultanovTzouTzou}, we derive 
	\begin{multline*}
		\int_{\mathbb{D}} L_a^{-1}\tilde{u}_j(s') \tilde{u}_j(s') e^{\tilde{\phi}(s')} ds' = \int_{\mathbb{D}} \left( u_j(x^*) + \varepsilon(b_1s_1 + b_2s_2) + R_2^2(\varepsilon s') \right) \times\\
		\times L_a^{-1} \left( u_j(x^*)e^{\phi(x^*)} + \varepsilon(c_1s_1 + c_2s_2) + R_2^1(\varepsilon s')\right) ds'
	\end{multline*}
	where $R_2^2$ is the reminder term of the Taylor series for $\tilde u_j$ near zero and $b_1$, $b_2$ are appropriate constants. Next, we note that 
	\begin{equation*}
		\int_{\mathbb{D}} L_a^{-1}[1](s')s_j ds' = 0, \qquad \|R^{j}_2(\varepsilon \cdot)\|_{H^1(\mathbb D)} = O(\varepsilon^2), 
	\end{equation*}
	for $j=1,$ $2$. Therefore, the penultimate identity gives
    \begin{multline*}
		\int_{\mathbb{D}} L_a^{-1}\tilde{u}_j(s') \tilde{u}_j(s') e^{\tilde{\phi}(s')} ds' = |u_j(x^*)|^2e^{\phi(x^*)} \int_{\mathbb D} L_a^{-1}[1](s')ds' + O(\varepsilon^2)\\ 
        = \frac{2\pi}{K_a}|u_j(x^*)|^2e^{\phi(x^*)} + O(\varepsilon^2).
    \end{multline*}

	Since $u_j$ is smooth on $\partial M$, it follows that $\tilde{u}_j(t) - \tilde{u}_j(0) = O_{H^{1/2}}(\varepsilon)$ as $\varepsilon\rightarrow 0$. Additionally, we recall that $R_{I,a}$, $R_{log,a}$, and $R_{\infty,a}$ are bounded operators from $H^{1/2}(\mathbb{D})$ to $H^{1/2}(\mathbb{D})^*$. Therefore, using (4.4) in \cite{NursultanovTzouTzou}, we write
	\begin{align*}
		\langle L_a^{-1} R_{log,a} L_a^{-1}\tilde{u}_j, \tilde{u}_je^{\tilde{\phi}} \rangle &=  \langle  R_{log,a} L_a^{-1}\tilde{u}_j, L_a^{-1}[\tilde{u}_je^{\tilde{\phi}}] \rangle\\
		&= |u_j(x^*)|^2e^{\phi(x^*)}  \langle  R_{log,a} L_a^{-1}[1], L_a^{-1}[1] \rangle + O(\varepsilon).
	\end{align*}
	Further, using \eqref{La u = 1}, we obtain
	\begin{align*}
		&\langle R_{log,a} L_a^{-1}\tilde{u}_j, L_a^{-1}[\tilde{u}_je^{\tilde{\phi}}] \rangle\\
		&= |u_j(x^*)|^2 e^{\phi(x^*)} \frac{a}{K_a^2} \int_{\mathbb{D}} \frac{1}{ (1-|s'|^2)^{1/2}} \int_\D  \frac{\log\left((t_1 - s_1)^2 + a^2 (t_2-s_2)^2 \right)^{1/2}}{ (1-|t'|^2)^{1/2}} dt' ds' + O(\varepsilon).
	\end{align*}
	Similarly, we collect expressions for $\langle L_a^{-1} R_{I,a} L_a^{-1}\tilde{u}_j, \tilde{u}_je^{\tilde{\phi}} \rangle $, $\langle R_{\infty,a} L_a^{-1}\tilde{u}_j, L_a^{-1}[\tilde{u}_je^{\tilde{\phi}}] \rangle$ and $\langle R_{F,a} L_a^{-1}\tilde{u}_j, L_a^{-1}[\tilde{u}_je^{\tilde{\phi}}] \rangle$ below
	\begin{multline*}
		\langle L_a^{-1} R_{I,a} L_a^{-1}\tilde{u}_j, \tilde{u}_je^{\tilde{\phi}} \rangle  = a |u_j(x^*)|^2 e^{\phi(x^*)} \int_{\mathbb{D}} L_a^{-1}[1](t') dt'  \int_{\mathbb{D}} L_a^{-1}[1](t') dt' + O(\varepsilon)\\
		=\frac{4\pi^2a}{K_a^2} |u_j(x^*)|^2 e^{\phi(x^*)} + O(\varepsilon),
	\end{multline*}
	\begin{multline*}
            \langle R_{\infty,a} L_a^{-1}\tilde{u}_j, L_a^{-1}[\tilde{u}_je^{\tilde{\phi}}] \rangle 
		= |u_j(x^*)|^2 e^{\phi(x^*)} \frac{a}{K_a^2} \times \\
            \times \int_{\mathbb{D}} \frac{1}{ (1-|s'|^2)^{1/2}} \int_\D \frac{(t_1 - s_1)^2 - a^2 (t_2 - s_2)^2}{(t_1 - s_1)^2 + a^2 (t_2 - s_2)^2} \frac{1}{ (1-|t'|^2)^{1/2}} dt' ds' + O(\varepsilon),
	\end{multline*}
	and finally (see page 10045 in \cite{NursultanovTradTzou}), 
	\begin{equation*}
		\langle R_{F,a} L_a^{-1}\tilde{u}_j, L_a^{-1}[\tilde{u}_je^{\tilde{\phi}}] \rangle = O(\varepsilon).
	\end{equation*}
	Using the identities above and \eqref{lambda minus lambda 1} we complete the proof.
\end{proof}

\section{Acknowledgement}
M.N. was partially supported by the grant of the Science Committee of the Ministry of Education and Science of the Republic of Kazakhstan, Grant No. AP14870361. M.N. was partially supported by the Academy of Finland, grants 353096 and 347715.

\bibliographystyle{plain}
\bibliography{references}

\setlength{\parskip}{0pt}





\end{document}